\documentclass[11pt]{article}
\usepackage{microtype}

\usepackage{amsmath,amssymb}
\usepackage{amsfonts,bm,bbm} 
\usepackage{amssymb,amsthm,enumitem}
\usepackage{tikz}
\usetikzlibrary{decorations.pathmorphing,positioning,patterns,snakes}
\usepackage{multirow}
\usepackage{xcolor,graphicx,float,booktabs}
\usepackage{tabularx}
\newcolumntype{L}[1]{>{\raggedright\arraybackslash}m{#1}}
\newcolumntype{C}[1]{>{\centering\arraybackslash}m{#1}}
\newcolumntype{R}[1]{>{\raggedleft\arraybackslash}m{#1}}

\usepackage{epstopdf}
\usepackage[algo2e,ruled	,vlined]{algorithm2e}

\usepackage{verbatim} 
\allowdisplaybreaks[4]
\usepackage{authblk}

\usepackage{hyperref}
\hypersetup{hypertex=true,
colorlinks=true,
linkcolor=blue,
anchorcolor=blue,
citecolor=blue}

\usepackage[round]{natbib}
\usepackage[
margin=1in,
includefoot,
footskip=30pt,
]{geometry}

\newcommand{\cF}{\mathcal{F}}

\newcommand{\cO}{\mathcal{O}}

\newcommand{\cX}{\mathcal{X}}
\newcommand{\E}{\mathbb{E}}
\newcommand{\R}{\mathbb{R}}
\newcommand{\Var}{\mathrm{Var}}

\newcommand{\bx}{\bf{x}}

\newtheorem{theorem}{Theorem}[section]
\newtheorem{proposition}[theorem]{Proposition}
\newtheorem{corollary}[theorem]{Corollary}
\newtheorem{example}[theorem]{Example}
\newtheorem{remark}{Remark}[section]

\newtheorem{assumption}{Assumption}[section]
\numberwithin{equation}{section}

\title{Enhanced Derivative-Free Optimization Using Adaptive Correlation-Induced Finite Difference Estimators}

\author[1]{Guo Liang\footnote{liangguo000221@ruc.edu.cn}}
\author[2]{Guangwu Liu\footnote{msgw.liu@cityu.edu.hk}}
\affil[1]{Institute of Statistics and Big Data\protect\\Renmin University of China\protect\\Beijing, China}

\author[1]{Kun Zhang\footnote{kunzhang@ruc.edu.cn}}

\affil[2]{Department of Management Sciences\protect\\City University of Hong Kong\protect\\Tat Chee Avenue, Kowloon, Hong Kong, China}

\begin{document}

\normalsize

\maketitle

\begin{abstract}

Gradient-based methods are well-suited for derivative-free optimization (DFO), where finite-difference (FD) estimates are commonly used as gradient surrogates. Traditional stochastic approximation methods, such as Kiefer-Wolfowitz (KW) and simultaneous perturbation stochastic approximation (SPSA), typically utilize only two samples per iteration, resulting in imprecise gradient estimates and necessitating diminishing step sizes for convergence. In this paper, we first explore an efficient FD estimate, referred to as correlation-induced FD estimate, which is a batch-based estimate. Then, we propose an adaptive sampling strategy that dynamically determines the batch size at each iteration. By combining these two components, we develop an algorithm designed to enhance DFO in terms of both gradient estimation efficiency and sample efficiency. Furthermore, we establish the consistency of our proposed algorithm and demonstrate that, despite using a batch of samples per iteration, it achieves the same convergence rate as the KW and SPSA methods. Additionally, we propose a novel stochastic line search technique to adaptively tune the step size in practice. Finally, comprehensive numerical experiments confirm the superior empirical performance of the proposed algorithm.

\end{abstract}

\section{Introduction}\label{sec:Introduction}

Stochastic optimization aims to find the minimization (or maximization) of a function in the presence of noise. Specifically, in this paper, we consider solving the following unconstrained stochastic optimization problem:
\begin{align}\label{eq:Problem}
	\min_{\bx \in \cX} F({\bx}) = \E[f({\bx})],
\end{align}
where $\cX \subseteq \R^{d}$ is a convex set, $F:\cX \rightarrow \R$ is the true performance, and $f$ is a noisy function. This problem has a wide range of application, including simulation optimization \citep[see, e.g.,][]{chang2013stochastic,hu2024convergence} and reinforcement learning \citep[see, e.g.,][]{fazel2018global,mania2018simple}. Among problem \eqref{eq:Problem}, a difficult but important case lies in the lack of the closed form of $F({\bx})$, and only estimates of the output function are available. That is, for any $\bx \in \cX$, we can only get an unbiased but noisy estimate of $F({\bx})$, i.e., $f({\bx})$. Such problem is derivative-free optimization (DFO, sometimes referred to as black-box optimization). As the problem becomes complex, the DFO will become increasingly important, and \cite{Golovin2017GoogleVizier} mention that ``any sufficiently complex system acts as a blackbox when it becomes easier to experiment with than to understand''.

Much literature has discussed the methodology development of DFO. The first category of algorithms are heuristic methods. For instance, the Nelder-Mead simplex algorithm, a direct-search-based method, is widely applied in practical scenarios \citep{barton1996nelder}. A key limitation of these algorithms is the lack of theoretical convergence guarantees \citep{spall2005introduction}. Another line involves transforming stochastic problems into deterministic ones, taking advantages of deterministic optimization. For examples, sample average approximation \citep{kim2015guide} generates many sample paths and uses the sample mean to estimate the unknown expectation; metamodels such as response surface methodology and Gaussian process are also used to fit the unknown function \citep{hong2021surrogate}. Recently, model-based trust region methods, which integrate model fitting and trust region techniques, have been developed to address simulation optimization problems. Notable examples include the stochastic trust-region response-surface method \citep{chang2013stochastic} and model-based trust-region derivative-free algorithms \citep{Shashaani2018}. Interested readers may refer to \cite{audet2017derivative} and \cite{larsonWild2019DFO} for a survey of these methods. 

Although our focus is on derivative-free approaches, \cite{audet2017derivative} comment that ``if gradient information is available, reliable and obtainable at reasonable cost, then gradient-based methods should be used.'' On the other hand, \cite{shi2023numerical} perform lots of experiments, showing the efficiency of gradient-based methods with finite-difference gradient surrogates. Therefore, we would like to consider gradient-based stochastic search algorithm to solve \eqref{eq:Problem} in this paper. Note that under some smoothness conditions (the assumption ensures the application of the gradient-based methods), solving \eqref{eq:Problem} is equivalent to finding the zero of the true gradient $\nabla F(\bx)$. Then, gradient-based methods are same as the stochastic approximation (SA) and take the general form:
\begin{align}\label{eq:SA}
	{\bx}_{k+1} = \Pi_{\cX}\left({\bx}_k - a_k \widehat{\nabla}F({\bx}_k)\right),
\end{align}
where ${\bx}_k$ and ${\bx}_{k+1}$ are the current solution and next prediction, respectively, $\widehat{\nabla}F({\bx}_k)$ is the estimate of the true gradient $\nabla F({\bx}_k)$ at $k$-th iteration, $a_k > 0$ is the step size, and $\Pi_{\cX}(\cdot)$ is the projection operator onto the feasible region $\cX$.

Under the black-box setting, the gradient estimate in \eqref{eq:SA} is usually substituted by finite difference (FD) gradient and in this case, \eqref{eq:SA} can date back to the Kiefer-Wolfowitz stochastic approximation (KWSA) algorithm \citep{Kiefer1952Stochastic}. Due to the variance of the FD estimate, the step size should tend to 0 to ensure the convergence. A similar method in high dimension is the simultaneous perturbation stochastic approximation (SPSA) method \citep{Spall1992Multivariate,Spall1997Multivariate}. In practice, the initial step size is crucial for the performance. For example, if the step size is too large, the solution may jump to a distant position, whereas if it is too small, the convergence rate of the KWSA and SPSA algorithms will be extremely slow, especially when the true performance is flat \citep[see Section 3.1.2 in][]{broadie2011general}.

To choose a proper step size, especially avoid the degradation of convergence rate, \cite{broadie2011general} present a scaled-and-shifted KW (SSKW) algorithm, adaptively tune the step size as well as the perturbation by introducing extra 9 hyperparameters. Firstly, they scale up the initial step size to avoid it is too small. Then, in the ``shifted'' procedure, they reduce the step size until the iteration falls into $\cX$. Another straightforward method lies in stochastic line search \citep[see, e.g.,][]{berahas2019derivative}, using the simulated function values to adjust the step size. Similarly, stochastic line search shrinks the step size until the stochastic Armijo condition is satisfied. Compared with SSKW method, stochastic line search is more flexible because the step size is not a fixed form, even not required to tend to 0. Consequently, this method necessitates the accuracy of the descent direction, and thus a batch of samples should be utilized at each iteration for gradient estimate. In this paper, we aim to introduce batch-based optimization to enhance the accuracy of gradient estimation, while using a larger step size to preserve the algorithm's effectiveness.

Batch-based optimization is inspired by the mini-batch method, which serves as an improvement over stochastic gradient descent or the Robbins-Monro (RM) algorithm \citep{Robbins1951Stochastic}. By utilizing a batch of samples to enhance the precision of stochastic gradient estimates, the mini-batch method is better suited for training neural networks. When using the FD gradient, the challenge lies in the construction of an accurate batch-based FD estimate. \cite{Fox1989Replication} and \cite{Zazanis1993Convergence} study the convergence of the standard batch-based FD estimator and provide the theoretically optimal perturbation size by minimizing the mean squared error (MSE). However, the optimal perturbation is related to the structure information of the blackbox and we do not know it when using the FD estimator. To overcome this issue, \cite{Li2020Optimally} propose a two-stage method, estimating the perturbation in the first stage and then generating remaining samples at the estimated perturbation to give the expectation-minimization FD (EM-FD) estimator. Recently, \cite{liang2024cor} propose a correlation-induced FD (Cor-FD) estimator, using all samples in a batch to estimate the perturbation and then recycling them to estimate the gradient. Cor-FD is available when the batch size is small and it is shown that Cor-FD possesses a reduced variance, and in some cases a reduced bias, compared to the optimal FD estimator.

Given that Cor-FD is efficient when the budget is limited, it is suitable for batch-based DFO algorithm \citep{Wang2024DFO}. However, selecting an appropriate batch size in each iteration is crucial. Specifically, if the batch size is too small, the descent direction may lack sufficient accuracy, limiting adjustments to minor corrections along this direction. Conversely, if the batch size is too large, samples may be wasted, as the descent direction does not require excessive precision. In fact, the most efficient algorithms should employ a progressive batching approach in which the batch size is initially small, and increases as the iteration progresses \citep{bollapragada2018progressive}. For this purpose, the adaptive sampling condition called the {\it{norm condition}} \citep[see, e.g.,][]{bollapragada2024derivative} has been proposed, which sets the batch size based on the {\it{signal-to-noise ratio}} (i.e., the ratio between the true gradient and the gradient estimate error).

By using the Cor-FD estimate and the adaptive sampling condition, the gradient estimate in \eqref{eq:SA} has been determined. Another challenge lies in the determination of the step size. The classical stochastic line search method introduces a relaxed Armijo condition \citep{berahas2019derivative,shi2023numerical}. However, step sizes that meet this condition are suboptimal and fall short of being ideal. Consequently, the algorithm converges only to a neighborhood of the optimal value, the size of which depends on the noise level in the black-box function \citep{berahas2019derivative}. Ensuring convergence requires allocating additional samples to determine a suitable step size. To address this, the paper proposes increasing the number of simulations for step sizes that meet the relaxed Armijo condition, enabling a more precise assessment of their suitability.

In summary, the main contributions of this paper are as follows:

\begin{enumerate}

	\item {\it{We propose an adaptive Cor-FD algorithm for DFO.}} We explore an efficient surrogate of $\widehat{\nabla}F({\bx}_k)$ in \eqref{eq:SA}: the Cor-FD estimate, and propose an adaptive sampling strategy for solving DFO. Specifically, the Cor-FD estimator exhibits a variance reduction property, enabling it to perform as well as, or even better than, the theoretically optimal FD estimate. That is, the Cor-FD estimator achieves a relatively small estimation error, and thus enhance the estimation of $\widehat{\nabla}F({\bx}_k)$. Furthermore, the proposed adaptive sampling strategy dynamically determines an appropriate batch size for constructing the Cor-FD estimate at each iteration, thereby further enhancing the efficiency and accuracy when solving DFO.
		
	\item {\it{We provide the consistency and convergence rate of our proposed algorithm.}} We prove that if the step size is set as some appropriate constants, the sequence of iterates generated by \eqref{eq:SA} converges linearly in expectation to the optimal solution. This result generalizes that of the deterministic (full-batch case in machine learning) gradient   descent. Moreover, it is surprising that although we consume a batch of samples to estimate the gradient, the algorithm maintains the optimal sample complexity in the sense that it reaches the optimal convergence rate of the KWSA algorithm.
	
	\item {\it{We revise the stochastic Armijo condition, establish an efficient DFO algorithm, and demonstrates its practical benefits through numerical experiments.}} We demonstrate that when step sizes are selected based on the relaxed Armijo condition, the algorithm converges only to a neighborhood of the optimal value. Subsequently, we heuristically propose new conditions to evaluate the suitability of step sizes. Using these conditions, we design an efficient algorithm and perform extensive numerical experiments. The results indicate that, compared to traditional methods, the new approach yields solutions closer to the optimal value. 

\end{enumerate}

The rest of the paper is organized as follows. In Section \ref{sec:BG}, we give some backgrounds about the gradient-based stochastic optimization and batch-based CFD estimate. Section \ref{sec:Alg} presents the adaptive sampling condition, the stochastic line search and the complete algorithms based on constant and line-search step sizes. In Section \ref{sec:Results}, we present the main results about the algorithm in Section \ref{sec:Alg}. Sections \ref{sec:Experiment} applies our algorithm to solve DFO problems, followed by conclusions in Section \ref{sec:Conclusion}. Proofs are provided in the appendix.

\section{Preliminaries}\label{sec:BG}
\subsection{Gradient-Based Stochastic Optimization}
To find the optimal solution $\bx^*$ of \eqref{eq:Problem}, the classical method is gradient-based stochastic search (also known as stochastic approximation). Specifically, $\bx^*$ can be obtained by the recursion \eqref{eq:SA}. Without loss of generality, we assume that $\cX = \R^{d}$ and remove the projection operator. Then, the recursion is
\begin{align}\label{eq:SA_R}
	{\bx}_{k+1} = {\bx}_k - a_k \widehat{\nabla}F({\bx}_k).
\end{align}

There are two pivotal roles in \eqref{eq:SA_R}: the step size $a_k$ and the gradient estimate $\widehat{\nabla}F({\bx}_k)$. If the unbiased gradient estimate $\nabla f({\bx}_k)$ can be obtained, then \eqref{eq:SA_R} is the RM algorithm, which is the origin of the stochastic gradient descent. In this paper, we assume that one can only get the zeroth-order information with noise and the first-order information is unavailable. In this case, the gradient estimate $\widehat{\nabla}F({\bx}_k)$ is usually substituted by the FD method. Such methods include the KW and SPSA methods, which obtain $\widehat{\nabla}F({\bx}_k)$ with only 2 samples at each iteration. Consequently, the variance of $\widehat{\nabla}F({\bx}_k)$ is large and the step size $a_k$ should tend to 0 to ensure the recursion goes towards the optimal solution $\bx^*$. Although it is shown that the convergence rate of the KW and SPSA algorithms can reach $\cO(1/k^{2/3})$, where $k$ is the iteration number (the sample complexity is similar because these 2 methods consume 2 samples at each iteration), the convergence rate may be unattainable in practice. Even when there is no noise, using the diminishing step size may lead to degeneration, as can be seen in Example \ref{exm:degenerate} \citep{broadie2011general}.

\begin{example}\label{exm:degenerate}
	Consider finding the infimum point of the deterministic function $f(x) = 0.001x^2$ with the KW method:
	\begin{align*}
		x_{k+1} = x_k - a_k \frac{f(x_k + h_k) - f(x_k - h_k)}{2h_k}.
	\end{align*}
	If we set $a_k = \theta_a/k$ and $h_k = \theta_h/k^{1/4}$ with $\theta_a = 1$ and $\theta_h = 1$, respectively, then the KW algorithm becomes $x_{k+1} = x_{k} (1 - 1/(500k))$. Starting with $x_1 = 1$, we have
	\begin{align*}
		x_{k+1} = \prod_{i=1}^{k}\left(1 - \frac{1}{500i}\right) = \exp\left(\sum_{i=1}^{k}\log\left(1 - \frac{1}{500i}\right)\right) \geq \exp\left(-\sum_{i=1}^{k}\frac{1}{500i}\right) \geq \cO\left(\frac{1}{k^{0.002}}\right),
	\end{align*}
	where the first inequality is because for any $x \in (0,1)$, $\log(1 - x) \geq -x$. 
	Therefore, the MSE cannot converge faster than $\cO\left(k^{-0.004}\right)$, which is significantly slower than the theoretically optimal rate of the KW algorithm $\cO\left(k^{-2/3}\right)$. The reason is that the step size is too small, which limits the algorithm's rate of descent.
\end{example}

To address the problem in Example \ref{exm:degenerate}, one approach is to carefully adjust the initial value of the step size, $\theta_a$, to prevent degeneration in the algorithm's convergence rate. Along this line, \cite{broadie2011general} propose the SSKW method, which adaptively tunes both the step size and perturbation by introducing 9 additional hyperparameters. Another approach is to use a constant step size, which demands a highly accurate gradient estimate, requiring more samples or a batch of samples to compute $\widehat{\nabla}F({\bx}_k)$. In this paper, we focus on the second method, where the first challenge lies in constructing the batch-based FD estimator.

\begin{remark}
	Estimating gradient with a batch of samples and then optimizing with the gradient-based method have been considered for a long time. In large-scale optimization, mini-batch method that utilizes a small batch of data to decrease the variance of stochastic gradient is thought as a promising alternative for traditional stochastic gradient descent. However, it is seldom used in the field of black-box optimization because of the bias inherent in the gradient estimate. Due to the linearity of the expectation, simply increasing the sample size often fails to reduce the bias in the gradient estimate, which results in the algorithm converging to the neighborhood of the optimal solution \citep{bollapragada2024derivative}. Therefore, the challenge of applying batch optimization in black-box environments lies in constructing an effective unbiased (at least asymptotically unbiased) gradient estimator.
\end{remark}

\subsection{Gradient Estimation}
In this section, we review the standard batch-based CFD estimator and then introduce a specific estimator -- the Cor-CFD estimator \citep{liang2024cor} which has the variance-reduction property and is suitable as the surrogate of $\widehat{\nabla}F({\bx}_k)$ in \eqref{eq:SA_R}. Throughout this section, we assume $d = 1$\footnote{When $d > 1$, we can estimate the gradient at each coordinate or choose only one coordinate as the descent component. Therefore, we set $d = 1$ for the save of simplicity.} and $x_0$ is the point of interest. 
\subsubsection{Standard Finite Difference}
To estimate the gradient $\nabla F(x_0)$, the CFD method chooses a perturbation $h$ and generates $n$ independent and identically distributed (i.i.d.) samples at $x_0 + h$ and $x_0 - h$, denoted by $\{f_1(x_0 + h),...,f_n(x_0 + h)\}$ and $\{f_1(x_0 - h),...,f_n(x_0 - h)\}$, respectively. Then, the CFD estimator is defined as
\begin{align*}
	g_{n,h}(x_0) = \frac{1}{n}\sum_{i=1}^{n}\frac{f_i(x_0 + h) - f_i(x_0 - h)}{2h}.
\end{align*}
If we assume that the function $F(x)$ is thrice continuously differentiable at $x_0$, then according to the Taylor expansion of $F(x_0 + h)$ and $F(x_0 - h)$ at $x_0$, we have
\begin{align}\label{eq:CFD_Bias}
	\E[g_{n,h}(x_0)] = \frac{F(x_0 + h) - F(x_0 - h)}{2h} = \nabla F(x_0) + B h^2 + o(h^2),
\end{align}
where $B = \nabla^3 F(x_0) /6$ is related to the third-order information. On the other hand, if we assume the noise of $f(x)$, denoted by $\sigma(x)$, is continuous at $x_0$, we have 
\begin{align}\label{eq:CFD_Variance}
	\Var[g_{n,h}(x_0)] = \frac{\sigma^2(x_0 + h) + \sigma^2(x_0 - h)}{4nh^2} = \frac{\sigma^2(x_0) + o(1)}{2nh^2}.
\end{align}

From \eqref{eq:CFD_Bias} and \eqref{eq:CFD_Variance}, it is obvious that $h$ controls the performance of the standard CFD estimator. When $h$ is too large, the bias will be unacceptable (see \eqref{eq:CFD_Bias}) and when $h$ is too small, the variance will explode (see \eqref{eq:CFD_Variance}). By minimizing the mean squared error (MSE, equal to $\mbox{Bias}^2 + \mbox{Variance}$), we find that the optimal perturbation is $h^* = \left(\sigma^2(x_0)/(4nB^2)\right)^{1/6}$.\footnote{Note that here we implicitly assume that $B \neq 0$. If $B = 0$, the CFD estimator is unbiased and the optimal perturbation is $\infty$. In this case, we prefer to choose a very large perturbation to approximate the optimal CFD estimator or utilize the forward and backward FD method.} However, we do not know the constants $B$ and $\sigma^2(x_0)$, leading the standard CFD method to be inefficient in practice. A common approach involves using pilot samples to estimate the unknown constants, constructing an estimate of the optimal perturbation, and then plugging it in the CFD estimator to generate the gradient estimator with new samples \citep{Li2020Optimally}. This method divides the total sample into two subsets: one for estimating the constants and the other for estimating the gradient. Consequently, when the sample size is limited, this method becomes inefficient because only a small portion of the sample is allocated for constant estimation, leading to significant errors in perturbation estimation. 

\subsubsection{Correlation-Induced Finite Difference}
The Cor-CFD method also includes two stages but utilizes the entire sample set for both perturbation estimation and gradient construction. Unlike previous approaches that divide the samples into separate subsets, this method leverages the same complete sample set for both tasks, thereby ensuring the accuracy of perturbation estimation while simultaneously maintaining the precision of the gradient estimation. 

Specifically, Cor-CFD method divides the total sample pairs $n$ by $K \times n_b$, where $K$ is the number of the perturbations generated from a user-specific distribution $\mathcal{P}$ (denoted by $h_1,...,h_K$), and $n_b$ is the number of the sample pairs at each perturbation (denoted by $\{f_1(x_0 + h_k), f_1(x_0 - h_k)\},...,\{f_{n_b}(x_0 + h_k), f_{n_b}(x_0 - h_k)\}$ at $h_k$, $k = 1,...,K$). Then, we have $K$ standard CFD estimators based on perturbations $h_1,...,h_K$, respectively, denoted by $g_{n_b,h_1}(x_0),..., g_{n_b,h_K}(x_0)$. From \eqref{eq:CFD_Bias} and \eqref{eq:CFD_Variance}, the expectation and variance of $g_{n,h}(x_0)$ are linear with respect to (w.r.t.) $(1, h^2)$ and $1/(nh^2)$, respectively. Therefore, using the bootstrap technique, we can get the estimates of $\E[g_{n_b,h_k}(x_0)]$ and $\Var[g_{n_b,h_k}(x_0)]$, denoted by $\E_*[g_{n_b,h_k}(x_0)]$ and $\Var_*[g_{n_b,h_k}(x_0)]$, respectively. Performing the linear regression gives the estimates $\left(\widehat{\widehat{\nabla}} F(x_0), \widehat{B}\right)$ and $\widehat{\sigma}^2(x_0)$. 

\begin{figure}[h!]
	\centering
	\caption{Comparison of the pilot samples, transformed samples and optimal samples for $f(x) = 10\sin(x) + \mbox{noise}$.}
	\hspace*{-2cm}
	\includegraphics[scale = 0.45]{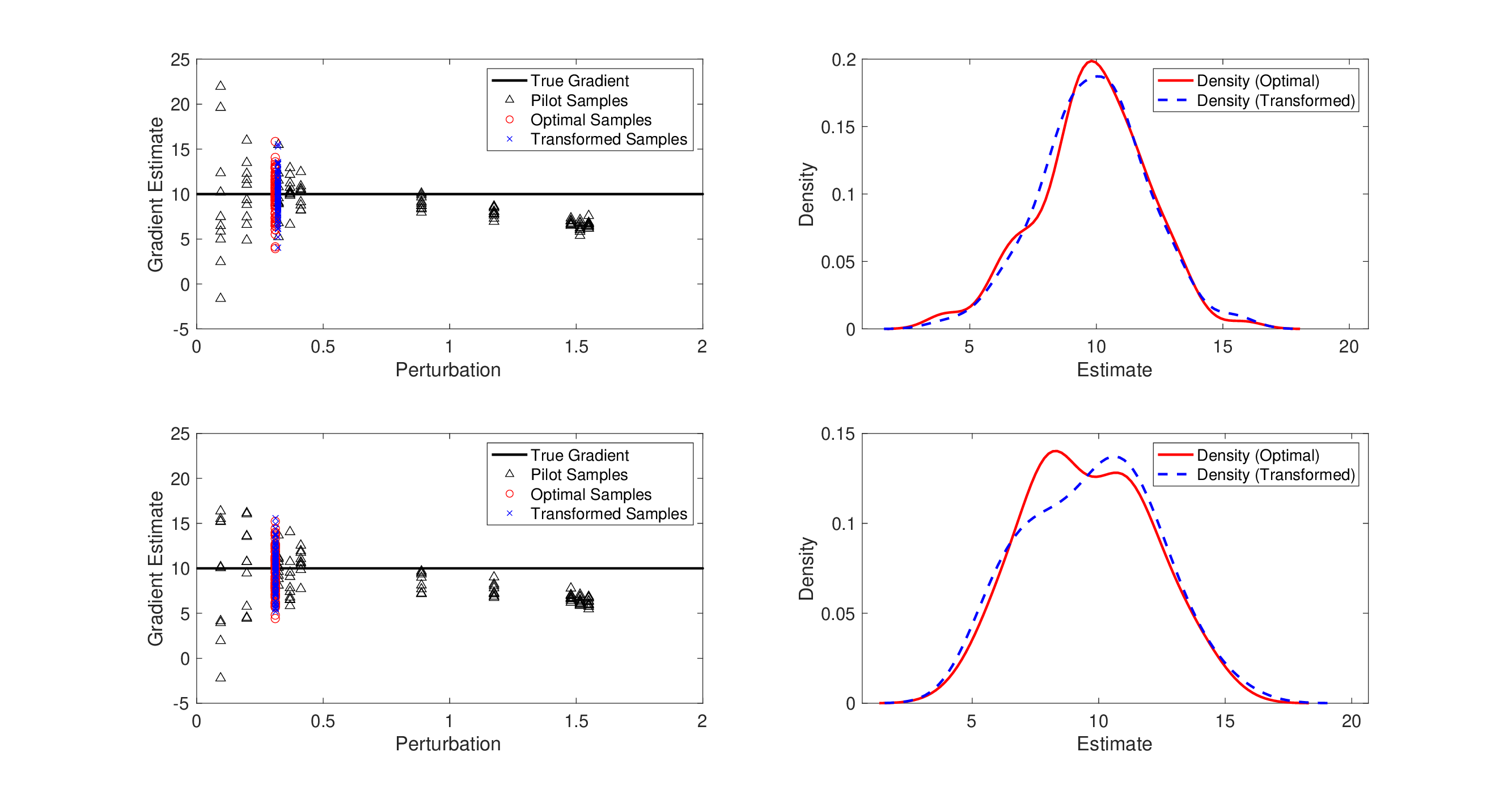}
	\vspace{-1cm}
	\label{fig:transform}
\end{figure}

Although the unknown constants $B$ and $\sigma^2(x_0)$ have been estimated, there is no extra sample to estimate the gradient, so the Cor-CFD method resorts to recycling the samples $\{f_1(x_0 + h_k), f_1(x_0 - h_k)\},...,\{f_{n_b}(x_0 + h_k), f_{n_b}(x_0 - h_k)\}$ for $k = 1,...,K$. Specifically, for a sample pair $\{f_j(x_0 + h_k), f_j(x_0 - h_k)\},\  j = 1,...,n_b, \ k = 1,...,K$, they use the transformation
\begin{align}\label{eq:transform}
	g_{h_k,j}^{cor}(x_0) = \frac{h_k}{\widehat{h}_n}\times\left(\frac{f_j(x_0 + h_k) - f_j(x_0 - h_k)}{2h_k} - \widehat{\widehat{\nabla}} F(x_0) - \widehat{B}h_k^2\right) + \widehat{\widehat{\nabla}} F(x_0) + \widehat{B}\widehat{h}_n^2
\end{align}
to get an {\it{adjusted version}} of the CFD estimate with a sample pair $\{f_j(x_0 + h_k), f(x_0 - h_k)\}$ and $\widehat{h}_n = \left(\widehat{\sigma}^2(x_0)/\left(4n \widehat{B}^2\right)\right)^{1/6}$. In \eqref{eq:transform}, $\widehat{\widehat{\nabla}} F(x_0) + \widehat{B}h_k^2$ and $\widehat{\widehat{\nabla}} F(x_0) + \widehat{B}\widehat{h}_n^2$ are used to adjust the expectation and $h_k/\widehat{h}_n$ is used to adjust the variance. This transformation is shown in Figure \ref{fig:transform}, in which we estimate the gradient of $f(x) = 10\sin(x) + \mbox{noise}$ at $x = 0$. 

The upper two subplots of Figure \ref{fig:transform} illustrate the case where the noise follows $\mathcal{N}(0,1)$, while the lower two subplots represent the case where the noise follows $\mathcal{U}(-\sqrt{3}, \sqrt{3})$. To estimate the gradient, we set $n=100$, $K=10$, and $n_b=10$. Here, `$\triangle$' represents the realistically generated CFD samples ($n_b$ samples at each perturbation), `$\mbox{o}$' indicates the CFD samples at the optimal perturbation ($n$ samples), and `$\text{x}$' denotes the samples obtained by transforming `$\triangle$' one by one using \eqref{eq:transform}. As shown in the left panel, `$\triangle$'s cannot be used as the gradient estimate directly due to their volatility or offset. However, they explicitly indicate the relationships like \eqref{eq:CFD_Bias} and \eqref{eq:CFD_Variance}, and can be used to estimate the optimal perturbation efficiently. Then, \eqref{eq:transform} transforms the `$\triangle$'s to `$\mbox{x}$'s and the transformed samples closely approximate the real samples, making them viable substitutes. The right panel displays the distributions of the optimal and transformed samples, estimated from $n$ data points. Although the two densities are not identical, they exhibit strong overall similarity. In the end, the Cor-CFD estimator is defined as the average of all {\it{adjusted version}} (i.e., `$\text{x}$' in Figure \ref{fig:transform}):
\begin{align*}
	g_n(x_0) = \frac{1}{n}\sum_{k=1}^{K}\sum_{j=1}^{n_b}g_{h_k,j}^{cor}(x_0). 
\end{align*}

\begin{proposition}[Theorem 4.1 in \cite{liang2024cor}]\label{pro:Cor-CFD}
	Assume that $F(x)$ is fifth differentiable at $x_0$ with non-zero fifth derivative, and the noise $\sigma(x) > 0$ is continuous at $x_0$. For any $k = 1,...,K$ $(K \geq 2)$, let $h_k = c_k n_b^{-1/10}$ $(c_k \neq 0)$ and for any $j \neq k$, $c_j \neq c_k$. If $n \to \infty$, then we have
	\begin{align}
	\E\left[g_n(x_0)\right] = \alpha'(x_0) + \left(\frac{B\sigma^2(x_{0})}{4n}\right)^{1/3} + \left(\frac{4B^2}{\sigma^2(x_0)}\right)^{1/6}\frac{D \Lambda}{\sqrt{K}}n^{-1/3} + o\left(n^{-1/3}\right),\label{eq:DSR_bias}\\
	\Var\left[g_n(x_0)\right] = \left(\frac{B^2 \sigma^4(x_0)}{2n^2}\right)^{1/3} + \left(\frac{B^2 \sigma^4(x_0)}{2}\right)^{1/3}\frac{q - K}{K}n^{-2/3} + o\left(n^{-2/3}\right),\label{eq:DSR_variance}
\end{align}
 where ${\boldsymbol{c}} = [|c_1|,...,|c_K|]^{\top}$, ${\boldsymbol{c}^4} = \left[c_1^4,...,c_K^4\right]^{\top}$, $\Lambda = \boldsymbol{c^{\top}Pc^4}$, ${\boldsymbol{P}}=\boldsymbol{I} - \boldsymbol{X}_e(\boldsymbol{X}_e^{\top}\boldsymbol{X}_e)^{-1}\boldsymbol{X}_e^{\top}$, the first and second columns of $\boldsymbol{X}_e$ are $[1,...,1]^{\top}$ and $[h_1^2,...,h_K^2]$, respectively, $q = \left\|{\rm{Diag}}(\boldsymbol{c^{-1}})\boldsymbol{Pc}\right\|_2^2$, ${\rm{Diag}}(\boldsymbol{c^{-1}}) = {\rm{Diag}}\left(1/|c_1|,...,1/|c_K| \right)$, and $||{\boldsymbol{v}}||_2 = \left(\sum_{k=1}^{K}v_k^2\right)^{1/2}$ for any ${\boldsymbol{v}} \in \R^K$.
\end{proposition}

From Proposition \ref{pro:Cor-CFD}, the second term on the right hand side (RHS) of \eqref{eq:DSR_bias} and the first term on the RHS of \eqref{eq:DSR_variance} are the same as the optimal bias and variance (this can be verified by taking $h^*$ back into \eqref{eq:CFD_Bias} and \eqref{eq:CFD_Variance}, respectively), and the asymptotic bias and variance of the Cor-CFD estimate are slightly different. Usually, the variance will be reduced because the projection matrix $\boldsymbol{P}$ will make $q$ in \eqref{eq:DSR_variance} smaller than $K$. On the other hand, it follows from \eqref{eq:DSR_bias} that the bias may also be reduced if the signs of $B$ and $D\Lambda$ are opposite. Therefore, the performance of Cor-CFD estimate is similar to (or even better than) the optimal case and is suitable for batch-based optimization.

\section{Proposed Algorithm}\label{sec:Alg}
Here and after, at $k$-th iteration, we denote $n_k$ by the number of sample pairs at each coordinate and $g_k({\bx}_k)$ by the corresponding Cor-CFD estimate. Back to \eqref{eq:SA}, although we have demonstrated an appropriate method to surrogate the gradient, there are still two questions should be addressed. The first question is how to set $n_k$ at $k$-th iteration and the second question is how to set the step size $a_k$. In Section \ref{sec:AS}, we consider the first question and propose an algorithm with constant step sizes, and in Section \ref{sec:LS}, we consider the second question and provide a more practical algorithm.
\subsection{Adaptive Sampling}\label{sec:AS}
For the first question, note that we aim to ensure that the estimated gradient aligns with the descending direction, meaning the angle between $g_k({\bx}_k)$ and $\nabla F({\bx}_k)$ is acute. However, this is not fully achievable due to the inherent uncertainty in gradient estimation. To address this, we incorporate the uncertainty by defining a confidence region for $\nabla F({\bx}_k)$, ensuring that all $d$-dimensional vectors within this region align with the descending direction. Specifically, let $\cF_k = \sigma\{{\bx}_1, {\bx}_2, ..., {\bx}_k\}$ be the $\sigma$-field generated by ${\bx}_1, {\bx}_2, ..., {\bx}_k$. Denote $\boldsymbol{b}_k = \E[g_k({\bx}_k) | \cF_k] - \nabla F({\bx}_k)$ and $\boldsymbol{\epsilon}_k = g_k({\bx}_k) - \E[g_k({\bx}_k) | \cF_k]$, which are the bias and noise terms, respectively. Consider the confidence region 
\begin{align*}
	[\boldsymbol{l},\boldsymbol{u}] := \left[g_k({\bx}_k) - \boldsymbol{b}_k - \sqrt{\E[\boldsymbol{\epsilon}_k\circ \boldsymbol{\epsilon}_k|\cF_k]}/\theta, g_k({\bx}_k) - \boldsymbol{b}_k + \sqrt{\E[\boldsymbol{\epsilon}_k\circ \boldsymbol{\epsilon}_k|\cF_k]}/\theta\right],
\end{align*}
where $\circ$ denotes the element-wise product, $\sqrt{\E[\boldsymbol{\epsilon}_k\circ \boldsymbol{\epsilon}_k|\cF_k]}$ is a $d$-dimensional vector and each element denotes the standard deviation of the corresponding element in $g_k({\bx}_k)$, and $\theta$ is a user-specified hyperparameter. Then, we let $\min\left\{\E[\boldsymbol{l}|\cF_k]^{\top}\nabla F({\bx}_k), \E[\boldsymbol{u}|\cF_k]^{\top}\nabla F({\bx}_k)\right\} \geq 0$ to ensure all $d$-dimensional vectors that lie in the confidence region are the decent direction. Without loss of generality, we assume all the elements of $\nabla F({\bx}_k)$ are positive. Therefore, we only require $\left(\nabla F({\bx}_k) - \sqrt{\E[\boldsymbol{\epsilon}_k\circ \boldsymbol{\epsilon}_k|\cF_k]}/\theta\right)^{\top}\nabla F({\bx}_k) \geq 0$, which can be derived from
\begin{align}\label{eq:AS_Condition}
	\E[||\boldsymbol{\epsilon}_k||^2 | \cF_k] \leq \theta^2 ||\nabla F({\bx}_k)||^2
\end{align}
using the Cauchy-Schwarz inequality. 

In fact, \eqref{eq:AS_Condition} represents an efficient sampling condition, indicating that the variance must be sufficiently small. This condition controls the {\it{noise-to-signal ratio}} in gradient estimation, thereby improving its reliability. \eqref{eq:AS_Condition} is called the {\it{norm condition}} and has been considered by \cite{bollapragada2018progressive} and \cite{bollapragada2024derivative}. To apply this condition, we need to identify surrogates for $\E[||\boldsymbol{\epsilon}_k||^2 | \cF_k]$ and $||\nabla F({\bx}_k)||$. Note that $\E[||\boldsymbol{\epsilon}_k||^2 | \cF_k]$ represents the sum of variances across all coordinates. we can employ sample variance to estimate the variance of each component of $g_k({\bx}_k)$ and subsequently $\E[||\boldsymbol{\epsilon}_k||^2 | \cF_k]$. Despite the existence of the correlation, it is efficient to use the sample variance because from Proposition \ref{pro:Cor-CFD}, the correlation tends to reduce the variance. For $||\nabla F({\bx}_k)||$, $||g_k({\bx}_k)||$ can be chosen as an appropriate surrogate. Specifically, the estimated version of \eqref{eq:AS_Condition} is 
\begin{align}\label{eq:AS_Condition_Esti}
	\frac{\sum_{i=1}^{d}\widehat{\sigma}_i^2}{n_k} \leq \theta^2 ||g_k({\bx}_k)||^2,
\end{align}
where $\widehat{\sigma}_i^2$ denotes the sample variance of the transformed samples \eqref{eq:transform} at $i$-th coordinate. Then the algorithm with constant step size is proposed (see Algorithm \ref{alg:cor_cfd_dfo_cs}).
Some key points of Algorithm \ref{alg:cor_cfd_dfo_cs} to note:
\begin{enumerate}
	\item The detailed description of the Cor-CFD method is omitted here, as the algorithm primarily focuses on optimization. The parameter settings of Cor-CFD method used in the experiments are provided in the appendix.
	\item From \cite{liang2024cor}, the sample pairs $n_k$ is an integer multiple of the number of pilot perturbations. In Algorithm \ref{alg:cor_cfd_dfo_cs}, if $\lfloor\sum_{i=1}^{d}\widehat{\sigma}_i^2/\left(\theta^2 ||g_k({\bx}_k)||^2\right)\rfloor + 1$ does not meet this condition,  it can be slightly increased to ensure compliance.
	\item Because $n_k$ is the sample pairs number at each coordinate, we should add $2dn_k$ when updating the number of function evaluations $s$.
	\item The parameters $a$ and $\theta$ should satisfy specific constraints to ensure the efficiency of the algorithm, which will be detailed in Section \ref{sec:Results}.
\end{enumerate}

\begin{algorithm2e}[b!]
\fontsize{10pt}{10pt}\selectfont
\caption{Cor-CFD-based DFO Algorithm with Constant Step Size}
\label{alg:cor_cfd_dfo_cs}
\KwIn{\\
\quad $\mathcal{S}$: Total number of function evaluations, \\
\quad $n_0$: Initial sample pairs, \\
\quad $s$: Number of function evaluations, \\
\quad ${\bf x}_0$: Starting point, \\
\quad $\theta$: Adaptive sampling threshold, \\
\quad $a$: Step length ($a > 0$).}
\KwOut{\\
\quad The ultimate estimate ${\bf x}_k$.}

\BlankLine
\textbf{Initialization:} Set $k \gets 0$, $s \gets 0$. \\

\While{$s < \mathcal{S}$}{
    \ForEach{coordinate $i = 1, ..., d$}{
        Compute gradient estimate $g_{k,i}({\bf x}_k)$ using Cor-CFD with $n_k$ sample pairs, where $g_{k,i}({\bf x}_k)$ denotes the $i$-th component of $g_k({\bf x}_k)$. \\
        Compute the sample variance $\widehat{\sigma}_i^2$.
    }
    \If{\eqref{eq:AS_Condition_Esti} does not hold}{
        Increase $n_k$ to $\lfloor\sum_{i=1}^{d}\widehat{\sigma}_i^2/\left(\theta^2 ||g_k({\bx}_k)||^2\right)\rfloor + 1$, where $\lfloor \cdot \rfloor$ represents the largest integer that does not exceed $\cdot$. \\
        Add $\lfloor\sum_{i=1}^{d}\widehat{\sigma}_i^2/\left(\theta^2 ||g_k({\bx}_k)||^2\right)\rfloor + 1 - n_k$ sample pairs to update gradient estimate $g_k({\bf x}_k)$. \\
        Set $n_k = \lfloor\sum_{i=1}^{d}\widehat{\sigma}_i^2/\left(\theta^2 ||g_k({\bx}_k)||^2\right)\rfloor + 1$.
    }
    Update ${\bf x}_{k+1} \gets {\bf x}_k - a g_k({\bf x}_k)$. \\
    Update $s \gets s + 2 d n_k$ and $k \gets k + 1$.
}
\Return ${\bf x}_k$
\end{algorithm2e}

\subsection{Stochastic Line Search}\label{sec:LS}
In practice, it is difficult to set an appropriate step size because the optimization problem is a blackbox. To address this issue, a stochastic line search method has been proposed to adjust the step size \citep{berahas2019derivative,shi2023numerical}, inspired by the backtracking line search in deterministic optimization. 

Specifically, the classical stochastic line search begins with an initial step size $a_k = \tilde{a}$ and determine whether 
\begin{align}\label{eq:ls_case1}
	f({\bx}_k - a_kg_k({\bx}_k)) > f({\bx}_k) - l_1 a_k ||g_k({\bx}_k)||^2 + 2 \sigma_f
\end{align}
holds, where $0 < l_1 < 1$ is a parameter. If \eqref{eq:ls_case1} holds, then the function value at the next predicted step is significantly larger than that at the current step and the step size should not be chosen. In this case, we shrink $a_k \to l_2 a_k$, where $l_1 < l_2 < 1$. The procedure will stop until \eqref{eq:ls_case1} does not hold. Note that $\sigma_f$ can be substituted by the upper bound of the function noise. For example, if it is very large, then \eqref{eq:ls_case1} will never hold and the stochastic line search outputs the constant step size $\tilde{a}$.

After the stochastic line search \eqref{eq:ls_case1}, we get a step size which is not too ``bad'' (a ``bad'' step size means that the next predicted value is much larger than the current value), but is not guaranteed to be ``good'' (a ``good'' step size means that the next predicted value is smaller than the current value). In other words, \eqref{eq:ls_case1} is conservative and only used to exclude the ``bad'' cases. Similar to \eqref{eq:ls_case1}, we present a new criterion:
\begin{align}\label{eq:ls_case2}
	f({\bx}_k - a_kg_k({\bx}_k)) \leq f({\bx}_k) - l_1 a_k ||g_k({\bx}_k)||^2 - 2\sigma_f.
\end{align}
The intuition of \eqref{eq:ls_case2} is that it indicates that the function value at the next predicted step is significantly smaller than that at the current step. 
Note that condition \eqref{eq:ls_case2} may not always hold, particularly when the current solution is near the optimal value. However, increasing simulation replications to reduce the black-box function's uncertainty will eventually satisfy condition \eqref{eq:ls_case2}. 

Specifically, there exists a value $N$ such that
\begin{align}\label{eq:ls_case2_Average}
	\frac{1}{N}\sum_{i=1}^{N}f_i({\bx}_k - a_kg_k({\bx}_k)) \leq \frac{1}{N}\sum_{i=1}^{N}f_i({\bx}_k) - l_1 a_k ||g_k({\bx}_k)||^2 - 2\frac{\sigma_f}{\sqrt{N}}
\end{align}
because when $N \to \infty$, it is identical to $F({\bx}_k - a_kg_k({\bx}_k)) \leq F({\bx}_k) - l_1 a_k ||g_k({\bx}_k)||^2$ which closes to the standard line search criterion and holds under mild conditions. Conditions \eqref{eq:ls_case2} and \eqref{eq:ls_case2_Average} require using many samples to assess the appropriateness of the step size. This effort is justified because, near the optimal solution, maintaining algorithmic progress demands highly accurate step-size selection. Note that when $N \to \infty$ and $a_k \to 0$, \eqref{eq:ls_case2_Average} always tends to hold. That is, if a step size is small enough, it is also a ``good'' step size (see Theorem \ref{thm:Converge_Iteration} for the theoretical evidences).
Thus, in practice, it is enough to evaluate condition \eqref{eq:ls_case2} for a finite number of times. If it remains unsatisfied, we shrink the step size. Additionally, we enforce a lower bound on the step size to avoid it becoming too small. The specific algorithm is shown in Algorithm \ref{alg:cor_cfd_dfo_ls}.

\begin{remark}
	Stochastic line search is similar to the ``shifted'' procedure in \cite{broadie2011general}, as both methods reduce a ``too large'' step size until it becomes appropriate. However, these two methods are fundamentally different. The ``shifted" procedure in \cite{broadie2011general} is based on KWSA, with the step size sequence defined as $\theta_a/k^{\gamma}$, where $\gamma > 0$ is a hyperparameter, and only $\theta_a$ is reduced. In contrast, stochastic line search is based on the standard line search used in deterministic optimization. In stochastic line search, the step size does not necessarily decrease at a fixed rate and it may decrease, not change or even increase. The only requirement is that the predicted value of the next iteration is smaller than the current function value to some extent. As a result, stochastic line search is more flexible.
\end{remark}

\begin{algorithm2e}[b!]
\fontsize{10pt}{10pt}\selectfont
\caption{Cor-CFD-based DFO Algorithm with Stochastic Line Search}
\label{alg:cor_cfd_dfo_ls}
\KwIn{\\
\quad $\mathcal{S}$: Total number of function evaluations, \\
\quad $n_0$: Initial sample pairs, \\
\quad $s$: Number of function evaluations, \\
\quad $(l_1,l_2,\underbar{a},N_0)$: Parameters of stochastic line search, \\
\quad ${\bf x}_0$: Starting point, \\
\quad $\theta$: Adaptive sampling threshold, \\
\quad $\tilde{a}$: Initial step size ($\tilde{a} > 0$).}
\KwOut{\\
\quad The ultimate estimate ${\bf x}_k$.}

\BlankLine
\textbf{Initialization:} Set $k \gets 0$, $s \gets 0$, $a_k \gets \tilde{a}$. \\

\While{$s < \mathcal{S}$}{
    \ForEach{coordinate $i = 1, ..., d$}{
        Compute gradient estimate $g_{k,i}({\bf x}_k)$ using Cor-CFD with $n_k$ sample pairs, where $g_{k,i}({\bf x}_k)$ denotes the $i$-th component of $g_k({\bf x}_k)$. \\
        Compute the sample variance $\widehat{\sigma}_i^2$.
    }
    \If{\eqref{eq:AS_Condition_Esti} does not hold}{
        Increase $n_k$ to $\lfloor\sum_{i=1}^{d}\widehat{\sigma}_i^2/\left(\theta^2 ||g_k({\bx}_k)||^2\right)\rfloor + 1$, where $\lfloor \cdot \rfloor$ represents the largest integer that does not exceed $\cdot$. \\
        Add $\lfloor\sum_{i=1}^{d}\widehat{\sigma}_i^2/\left(\theta^2 ||g_k({\bx}_k)||^2\right)\rfloor + 1 - n_k$ sample pairs to update gradient estimate $g_k({\bf x}_k)$. \\
        Set $n_k = \lfloor\sum_{i=1}^{d}\widehat{\sigma}_i^2/\left(\theta^2 ||g_k({\bx}_k)||^2\right)\rfloor + 1$.
    }
    \While{\eqref{eq:ls_case1} holds}{
    	Shrink $a_k \gets l_2 a_k$.
    }
    \While{$a_k > \underbar{a}$}{
    	\If{\eqref{eq:ls_case2_Average} holds for some $N \leq N_0$}{
    		Break the loop.
    	}
    	\Else{
    		Shrink $a_k \gets l_2 a_k$.
    	}
    }
    Count the number of function evaluations for stochastic line search $n_{ls}$. \\
    Update ${\bf x}_{k+1} \gets {\bf x}_k - a_k g_k({\bf x}_k)$. \\
    Update $s \gets s + 2 d n_k + n_{ls}$ and $k \gets k + 1$.
}
\Return ${\bf x}_k$
\end{algorithm2e}

\section{Theoretical Results}\label{sec:Results}
In this section, we present the convergence results of Algorithm \ref{alg:cor_cfd_dfo_cs} and stochastic line search \eqref{eq:ls_case1}. Firstly, we state some assumptions, which suppose that the objective function $F(\bx)$ is smooth and strongly convex, satisfying the following regularity conditions.
\begin{assumption}[Differentiability]\label{ass:Diff}
	The function $F({\bx})$ is fifth continuously differentiable on $\cX$ and $\nabla^5 F({\bx}) \neq 0$.
\end{assumption}

\begin{assumption}[Lipschitz smoothness]\label{ass:LS}
	The gradient of $F({\bx})$ is $M$-Lipschitz, i.e., there exists a constant $M > 0$ such that $||\nabla{F}({\bx_1}) - \nabla{F}({\bx_2})|| \leq M ||{\bx_1} - {\bx_2}||$ for all ${\bx_1}, {\bx_2} \in \cX$, where $|| \cdot ||$ denotes the Euclidean norm.
\end{assumption}

\begin{assumption}[Strongly convex]\label{ass:SC}
	There exists a constant $m > 0$ such that $\lambda({\bx}) \geq m$ for all ${\bx} \in \cX$, where $\lambda({\bx})$ denotes the smallest eigenvalue of the Hessian matrix $H({\bx}) := \nabla^2 F({\bx})$.
\end{assumption}

Assumption \ref{ass:Diff} is a relatively strong condition, which is used for the theoretical completeness of the Cor-CFD method. Note that when $\cX$ is a compact set, Assumption \ref{ass:LS} holds automatically if Assumption \ref{ass:Diff} holds. Assumption \ref{ass:SC} ensures that \eqref{eq:Problem} has a unique solution ${\bx}^*$. Assumptions \ref{ass:LS} and \ref{ass:SC} are standard conditions when studying the convergence results of the gradient-based method \citep[see, e.g.,][]{ghadimi2012optimal,scheinberg2022finite}. Specifically, these two assumptions mean that $M$ and $m$ are the upper and lower bounds for all eigenvalues of the Hessian matrix.

In the following, we present Theorem \ref{thm:Converge_Iteration} to show the convergence result of Algorithm \ref{alg:cor_cfd_dfo_cs} and the proof is provided in Appendix \ref{app:proofTheorem2}. Note that in Theorem \ref{thm:Converge_Iteration}, we use the condition $\max\{\E[||\boldsymbol{b}_k||^2 | \cF_k], \E[||\boldsymbol{\epsilon}_k||^2 | \cF_k]\} \leq \theta^2 ||\nabla F({\bx}_k)||^2$ which is different from \eqref{eq:AS_Condition}. This condition is introduced solely for the convenience of the proof. It is reasonable because if we use Cor-CFD gradient estimate, $\E[||\boldsymbol{b}_k||^2 | \cF_k]$ and $\E[||\boldsymbol{\epsilon}_k||^2 | \cF_k]$ have the same order (see Proposition \ref{pro:Cor-CFD}). 
\begin{theorem}\label{thm:Converge_Iteration}
	Suppose that Assumptions \ref{ass:Diff}, \ref{ass:LS} and \ref{ass:SC} hold. Let ${\bx}_0$ be the initial point and at $k$-th iteration, let $\max\{\E[||\boldsymbol{b}_k||^2 | \cF_k], \E[||\boldsymbol{\epsilon}_k||^2 | \cF_k]\} \leq \theta^2 ||\nabla F({\bx}_k)||^2$, where $0 < \theta < m/(2M)$ is a threshold. If $0 < a_k = a \leq 1 / ((2\theta^2 + 2\theta + 1) M)$ for any $k \geq 0$, then we have
	\begin{align}\label{eq:Linear_Converge}
		\E\left[||{\bx}_k - {\bx}^*||^2\right] \leq \left(1 - (m - 2\theta M)a\right)^{k} ||{\bx}_0 - {\bx}^*||^2. 
	\end{align}
\end{theorem}

Observe that $0 < (m - 2\theta M)a < 1$ when $0 < \theta < m/(2M)$ and $a \leq 1 / ((2\theta^2 + 2\theta + 1) M)$. Under these conditions, ${\bx}_k$ converges linearly in expectation to the minimum point ${\bx}^*$. Theorem \ref{thm:Converge_Iteration} generalizes the convergence result of standard gradient descent. Notably, as $\theta \to 0$, the gradient estimate in each iteration approaches the true gradient at the current solution. Consequently, the RHS on \eqref{eq:Linear_Converge} converges to $(1 - ma)^{k} ||{\bx}_0 - {\bx}^*||^2$, where $0 < a < 1/M$. The optimal case occurs as $a \to 1/M$, with $\E[||{\bx}_k - {\bx}^*||^2]$ converging to 0 at a rate comparable to a geometric series, featuring an exponent approaching $m/M$. This result aligns perfectly with deterministic gradient descent. While Theorem \ref{thm:Converge_Iteration} assumes $0 < \theta < m/(2M)$, this condition can be relaxed to $0 < \theta < m/M$. Under this relaxation, a larger gradient estimation error is permissible, but a smaller step size is required, leading to a slower convergence rate (see Corollary \ref{cor:Converge_Iteration}).


\begin{corollary}[A weak version of Theorem \ref{thm:Converge_Iteration}]\label{cor:Converge_Iteration}
	Suppose that Assumptions \ref{ass:Diff}, \ref{ass:LS} and \ref{ass:SC} hold. Let ${\bx}_0$ be the initial point and at $k$-th iteration, let $\max\{\E[||\boldsymbol{b}_k||^2 | \cF_k], \E[||\boldsymbol{\epsilon}_k||^2 | \cF_k]\} \leq \theta^2 ||\nabla F({\bx}_k)||^2$, where $0 < \theta < m/M$ is a threshold. If $0 < a_k = a \leq (m - \theta M) / ((2\theta^2 + 2\theta + 1) M^2)$ for any $k \geq 0$, then we have
	\begin{align}\label{eq:Linear_Converge_Coro}
		\E\left[||{\bx}_k - {\bx}^*||^2\right] \leq (1 - (m - \theta M) a)^{k} ||{\bx}_0 - {\bx}^*||^2. 
	\end{align}
\end{corollary}

The proof of Corollary \ref{cor:Converge_Iteration} is provided in Appendix \ref{app:proofCorollary}. According to Corollary \ref{cor:Converge_Iteration}, the optimal scenario occurs when the constant step size is given by $a = (m - \theta M) / ((2\theta^2 + 2\theta + 1) M^2)$. While Corollary \ref{cor:Converge_Iteration} guarantees that ${\bx}_k$ converges linearly in expectation to the minimum point ${\bx}^*$, the convergence rate is slower compared to Theorem \ref{thm:Converge_Iteration}, as the step size is smaller. However, the step size proposed in the corollary is well-suited for larger thresholds, reducing the number of samples needed to meet the adaptive sampling condition. 

It is important to note that focusing solely on the convergence rate of ${\bx}_k$ is insufficient. For instance, when $\theta \to 0$, the gradient descent result is recovered (at least, the convergence rate is faster), but this requires the batch size to be infinity at each iteration to satisfy the adaptive sampling condition.
Therefore, it is crucial to consider both the convergence rate of ${\bx}_k$ and the associated sample complexity, i.e., the total stochastic function evaluations required to get an $\epsilon$-accurate solution. We employ the metric that ${\bx}_k$ is said to be an $\epsilon$-accurate solution if $\E\left[||{\bx}_k - {\bx}^*||^2\right] \leq \epsilon$. Note that the number of stochastic function evaluations at any iteration $k$ is $\mathcal{S}_k = 2dn_k$, where $n_k$ is the number of sample pairs at each coordinate. In the following, we present Theorem \ref{thm:Sample_Complex} to analyze the sample complexity of Algorithm \ref{alg:cor_cfd_dfo_cs} and the proof is provided in Appendix \ref{app:proofTheorem3}.

\begin{theorem}[Sample complexity]\label{thm:Sample_Complex}
	Under the same conditions as those in Theorem \ref{thm:Converge_Iteration}. Let $d = \cO(1)$. Denote $\mathcal{S}(\epsilon)$ by the total stochastic function evaluations to get an $\epsilon$-accurate solution. If the third derivative $\nabla^3 F({\bx})$ and the function noise $\sigma({\bx})$ are bounded, then we have
	\begin{align*}
		\mathcal{S}(\epsilon) \geq \mathcal{C}_1 \epsilon^{-3/2} + \mathcal{C}_2,
	\end{align*}
	where $\mathcal{C}_1$ and $\mathcal{C}_2$ are constants that depends on the threshold $\theta$, step size $a$, problem dimension $d$, unknown function and the simulation error. Specifically,
	\begin{align*}
		\mathcal{C}_1 = \frac{4d(\bar{\mathcal{C}}M^2)^{-3/2}}{3\bar{a}^{-3/2}\log(1/\bar{a})},\quad \mathcal{C}_2 = \frac{4d(\bar{\mathcal{C}}M^2||{\bx}_0 - {\bx}^*||^2)^{-3/2}}{3\log(1/\bar{a})},
	\end{align*}
	where $\bar{a} = 1 - (m - 2\theta M)a$, $\bar{\mathcal{C}} = \theta^2/\E[\mathcal{C}]$ and $\max\{\E[||\boldsymbol{b}_k||^2 | \cF_k], \E[||\boldsymbol{\epsilon}_k||^2 | \cF_k]\} \leq \mathcal{C}n_k^{-2/3}$.
\end{theorem}

According to Theorem \ref{thm:Sample_Complex}, achieving an $\epsilon$-accurate solution requires at least $\cO\left(\epsilon^{-3/2}\right)$ function evaluations. In other words, for a given total budget $\mathcal{S}$, the MSE of our algorithm is $\cO\left(\mathcal{S}^{-2/3}\right)$\footnote{Because $d$ is assumed to be $\cO(1)$, it only appears on constants $\mathcal{C}_1$ and $\mathcal{C}_2$, and has little influence on the convergence rate. When $d$ is related to the total sample budget, the sample complexity is different from that in Theorem \ref{thm:Sample_Complex} and we do not consider it in this paper.}, matching the optimal performance of the KW algorithm \citep{hu2024convergence}. This result stems from combining the Cor-CFD estimate with the adaptive sampling condition, which ensures that an appropriate number of samples are generated at each iteration, maximizing the use of sample information. For a fixed total budget, reliable gradient estimation allows the algorithm to maintain consistently sufficient descent, even with a reduced number of iterations. Furthermore, because the step size does not need to approach 0, the algorithm circumvents the degeneration scenario illustrated in Example \ref{exm:degenerate}.

In the following, we present the convergence result when using \eqref{eq:ls_case1} in the optimization algorithm. Our analysis relies on the following extra assumption.
\begin{assumption}[Boundedness of function noise]\label{ass:Bound_Error_Fun}
	There is a constant $\bar{\epsilon}_f > 0$ such that $|f({\bx}) - F({\bx})| \leq \bar{\epsilon}_f$ for all ${\bx} \in \cX$.
\end{assumption}
Although Assumption \ref{ass:Bound_Error_Fun} does not hold exactly in general, it can be ensured to hold with a high probability. For instance, if the noise in the function $f({\bx})$ follows a normal distribution $\mathcal{N}(F({\bx}),\sigma^2({\bx}))$ and the upper bound of the standard deviation $\sigma({\bx})$ is denoted by $\sigma$. Then, setting $\bar{\epsilon}_f = 3\sigma$ ensures that Assumption \ref{ass:Bound_Error_Fun} holds with a probability exceeding 99\%. Moreover, it is easy to observe that if condition \eqref{eq:ls_case2} holds, the current step size is guaranteed to decrease the objective function under Assumption \ref{ass:Bound_Error_Fun}.
Next, we present the line-search result with condition \eqref{eq:ls_case1}, i.e., Theorem \ref{thm:Converge_LS}. In Theorem \ref{thm:Converge_LS}, we assume that $||g_k({\bx}_k) - \nabla F({\bx}_k)|| \leq \theta |\nabla F({\bx}_k)|$, which is stronger than the adaptive sampling condition that appears in Theorem \ref{thm:Converge_Iteration}.\footnote{To obtain a similar result to Theorem \ref{thm:Converge_LS}, it is enough to assume that the gradient estimate error is bounded like \cite{berahas2019derivative} and the boundedness condition can be ensured to hold with a high probability by using the adaptive sampling condition \eqref{eq:AS_Condition} and Markov's inequality.} 

\begin{theorem}\label{thm:Converge_LS}
	Suppose that Assumptions \ref{ass:Diff}, \ref{ass:LS}, \ref{ass:SC} and \ref{ass:Bound_Error_Fun} hold. Let ${\bx}_0$ be the initial point and at $k$-th iteration, $a_k$ is chosen by \eqref{eq:ls_case1} from $\tilde{a}$. Let $||g_k({\bx}_k) - \nabla F({\bx}_k)|| \leq \theta |\nabla F({\bx}_k)|$, where $0 < \theta < 1$ is a threshold. If $0 < l_1 \leq (1 - \theta) / (2 (1 + \theta)) \leq l_2 < 1$, then we have $a_k > l_2/M$ and
	\begin{align}\label{eq:LS_Results}
		F({\bx}_k) - F({\bx}^*) - \frac{4\bar{\epsilon}_f}{1 - \rho} \leq \rho^{k} \left(F({\bx}_0) - F({\bx}^*) - \frac{4\bar{\epsilon}_f}{1 - \rho}\right),
	\end{align}
	where $\rho = 1 - \min\{2ml_1\tilde{a}(1-\theta)^2, 2ml_1l_2(1 - \theta)^2/M\}$.
\end{theorem}

The proof of Theorem \ref{thm:Converge_LS} is provided in Appendix \ref{app:proofTheorem4}. From Theorem \ref{thm:Converge_LS}, under specific assumptions, the classical stochastic line search method ensures the existence of a step size larger than $l_2/M$. Consequently, the line search procedure can terminate within a finite number of steps. However, the step size obtained through this process does not guarantee a monotonic decrease in the function value. Instead, the function value tends to fluctuate and converge within a neighborhood of the optimal value, rather than achieving precise convergence to the optimum (see \eqref{eq:LS_Results}). This limitation arises because the optimization problem involves a noisy black-box function. In such cases, the presence of inherent noise hinders the accurate identification of the optimal step size within a limited number of iterations, leading to inevitable fluctuations.

\section{Numerical Experiments}\label{sec:Experiment}
In this section, we test the performance of our algorithm on some examples. It is important to note that the examples presented may not strictly satisfy the assumptions provided in Section \ref{sec:Results}. For instance, if the support set of the black-box noise distribution is unbounded, Assumption \ref{ass:Bound_Error_Fun} does not hold. However, the assumptions we present are primarily intended to support the theoretical results. Even when these assumptions are violated, the algorithm still performs well under certain ``stress tests''. In this section, the main focus is on the algorithm's output, and we are less concerned with whether the tested function satisfies all the assumptions exactly. 

\subsection{Test Problems}
We consider three types of problems which are all constructed from deterministic optimization with added noise. We choose these test problems because when comparing with real simulation experiments, the true values of our test problem can be obtained easily, allowing us to investigate the strengths and limitations of the algorithm \citep{chang2013stochastic}. The first type of problem is a simple power function $F(x) = x^4$, but we add varying levels of noise to control the difficulty of the optimization problem. Specifically, we denote $f(x) = F(x) + \mathcal{N}(0,\sigma(x)^2)$, where $\sigma(x) = 0.1, 1, 10$. This problem is also used as a test for the performance of the SSKW algorithm \citep{broadie2011general}.
The second type of problem is the Rosenbrock function defined as
\begin{align*}
	F(x_1,x_2) = 100(x_2 - x_1^2)^2 + (x_1 - 1)^2.
\end{align*}
The unique minimum $(1, 1)$ lies in a narrow, parabolic valley $(x_1, x_1^2)$. Although the valley is easy to find, it is  difficult to find the minimum even when there is no noise. In this problem, we set $f(x_1, x_2) = F(x_1, x_2) + \mathcal{N}(0,1)$. The last type of problem is a 64-dimensional function \citep{schittkowski2012more} defined as
\begin{align*}
	F({\bx}) = \sum_{i=1}^{32} \left[10(x_{2i} - x_{2i-1})^2 + (1 - x_{2i-1})^2 \right]^4.
\end{align*}
This function is steep when far from the optimal solution and becomes flat near the optimum. In this problem, we set $f({\bx}) = F({\bx}) + \mathcal{N}(0,\sigma({\bx})^2)$, where $\sigma({\bx}) = 0.1, 1, 10$. 

\subsection{Performance Measures}
For the first type of problem ($d = 1$), we compare Algorithm \ref{alg:cor_cfd_dfo_ls} with KWSA algorithm. For other types of problem, we compare Algorithm \ref{alg:cor_cfd_dfo_ls} with SPSA algorithm. We present all the involved performance measures in the following.
\begin{itemize}
	\item {\it{Solution error.}} The solution error is defined as $||{\bx}_k^* - {\bx}^*||$, where ${\bx}^*$ and ${\bx}_k^*$ are the optimal solution and the last solution before the algorithm terminates, respectively. 
	\item {\it{Optimality gap.}} The optimality gap (OG) is defined as $F({\bx}_k^*) - F({\bx}^*)$. 
	\item {\it{Length of oscillatory period.}} The oscillatory period is defined as the number of sample pairs until the algorithm stops oscillating between boundary points. Mathematically, it is represented by the cardinality of the set $ \{k \geq 2 : {\bx}_{k-1},{\bx}_k \in \partial {\cX} \ \text{and}\  {\bx}_k \neq {\bx}_{k-1}\}$, where $\partial {\cX}$ is the boundary of $\cX$.
\end{itemize}

\subsection{Parameter Settings}
The parameter settings are as follows:
\begin{itemize}
	\item {\it{Parameters in Algorithm \ref{alg:cor_cfd_dfo_ls}.}} For all problems, we set the initial sample pairs $n_0 = 10$. The parameters of stochastic line search are set as $l_1 = 10^{-4}$, $l_2 = 0.5$, $\underbar{a} = 0$ and $N_0 = 10$. The adaptive sampling threshold is set as $\theta = 0.7$. The initial step size is set as $\tilde{a} = 1$.
	\item {\it{Parameters in KWSA.}} When using KWSA, we set the step size and perturbation as $\theta_a/k$ and $\theta_c/k^{1/4}$, respectively, which are suggested in \cite{broadie2011general}. The results of $\theta_a = \theta_c = 1$ are reported.
	\item {\it{Parameters in SPSA.}} The SPSA algorithm is instructed by \url{http://www.jhuapl.edu/SPSA/Pages/MATLAB.htm}. Specifically, we set the step size and perturbation as $\theta_a/(k + 50)^{0.602}$ and $\theta_c/k^{0.101}$. We run 20 macroreplications for $\theta_a = \{10^{-9}, 10^{-8}, ..., 1, 10, 100\}$, $\{\theta_c = 10^{-4}, 10^{-3}, ..., 1, 10, 100\}$ with $2000d$ function evaluations and obtain the parameters $(\theta_a, \theta_c)$ with minimal average OG. 
\end{itemize}

\subsection{Results}
Algorithm \ref{alg:cor_cfd_dfo_ls} is marked as ``AdaDFO'' across all the experiments.
For the first type of problem, we set $\cX = [-50, 50]$, $x_0 = 30$, and the results are reported in Table \ref{tab:Power4}. As shown in the table, the solution error of our algorithm decreases with an increasing function evaluation budget. It is important to note that when $\sigma(x) = 10$, the solution error does not decrease with an increasing sample size. This is because as the noise level rises, the sample size needed to satisfy the adaptive sampling condition also increases (see \eqref{eq:AS_Condition_Esti}). As the sample pairs grow from 100 to 1,000, the adaptive sampling condition may be not satisfied, causing the solution to fluctuate within a small range rather than strictly decrease. Furthermore, the table reveals that the KWSA algorithm oscillates at the boundaries and begins to converge after about 5000 iterations. This suggests that the traditional KWSA algorithm requires careful parameters tuning before it can be effectively applied. In contrast, the length of  oscillatory period of our algorithm remains zero, indicating steady descent and demonstrating the algorithm's reliability.

\begin{table}[t!]
\renewcommand\arraystretch{0.9}
\centering
\caption{Comparison of Algorithm \ref{alg:cor_cfd_dfo_ls} with KWSA algorithm for first type of problem.}
\label{tab:Power4}
\vspace{1em}
\begin{tabular}{cccccccc}
\toprule

 &  & \multicolumn{3}{c}{Solution error} & \multicolumn{3}{c}{Length of oscillatory period} \\
\cmidrule(lr){3-5} \cmidrule(lr){6-8}
$\sigma(x)$ & Method & \multicolumn{1}{r}{100} & \multicolumn{1}{r}{1,000} & \multicolumn{1}{r}{10,000} & \multicolumn{1}{r}{5\%} & \multicolumn{1}{r}{Median} & \multicolumn{1}{r}{95\%} \\
\midrule
\multirow{2}{*}{0.1} & AdaDFO & \multicolumn{1}{r}{0.18} & \multicolumn{1}{r}{0.12} & \multicolumn{1}{r}{0.10} & \multicolumn{1}{r}{0} & \multicolumn{1}{r}{0} & \multicolumn{1}{r}{0} \\
& KWSA       & \multicolumn{1}{r}{50.00}  & \multicolumn{1}{r}{50.00}  & \multicolumn{1}{r}{0.42} & \multicolumn{1}{r}{5000} & \multicolumn{1}{r}{5000} & \multicolumn{1}{r}{5000} \\
\multirow{2}{*}{1}   & AdaDFO & \multicolumn{1}{r}{0.23} & \multicolumn{1}{r}{0.20} & \multicolumn{1}{r}{0.14} & \multicolumn{1}{r}{0} & \multicolumn{1}{r}{0} & \multicolumn{1}{r}{0} \\
& KWSA       & \multicolumn{1}{r}{50.00}  & \multicolumn{1}{r}{50.00}  & \multicolumn{1}{r}{0.42} & \multicolumn{1}{r}{4999} & \multicolumn{1}{r}{5000} & \multicolumn{1}{r}{5000} \\
\multirow{2}{*}{10}  & AdaDFO & \multicolumn{1}{r}{0.35} & \multicolumn{1}{r}{0.38} & \multicolumn{1}{r}{0.33} & \multicolumn{1}{r}{0} & \multicolumn{1}{r}{0} & \multicolumn{1}{r}{0} \\
& KWSA       & \multicolumn{1}{r}{50.00}  & \multicolumn{1}{r}{50.00}  & \multicolumn{1}{r}{0.38} & \multicolumn{1}{r}{4999} & \multicolumn{1}{r}{4999} & \multicolumn{1}{r}{5000} \\

\bottomrule
\end{tabular}
\end{table}

\begin{figure}[b!]
	\centering
	\caption{Comparison of the average OG of Algorithm \ref{alg:cor_cfd_dfo_ls} and SPSA for second type of problem.}
	\includegraphics[scale = 0.75]{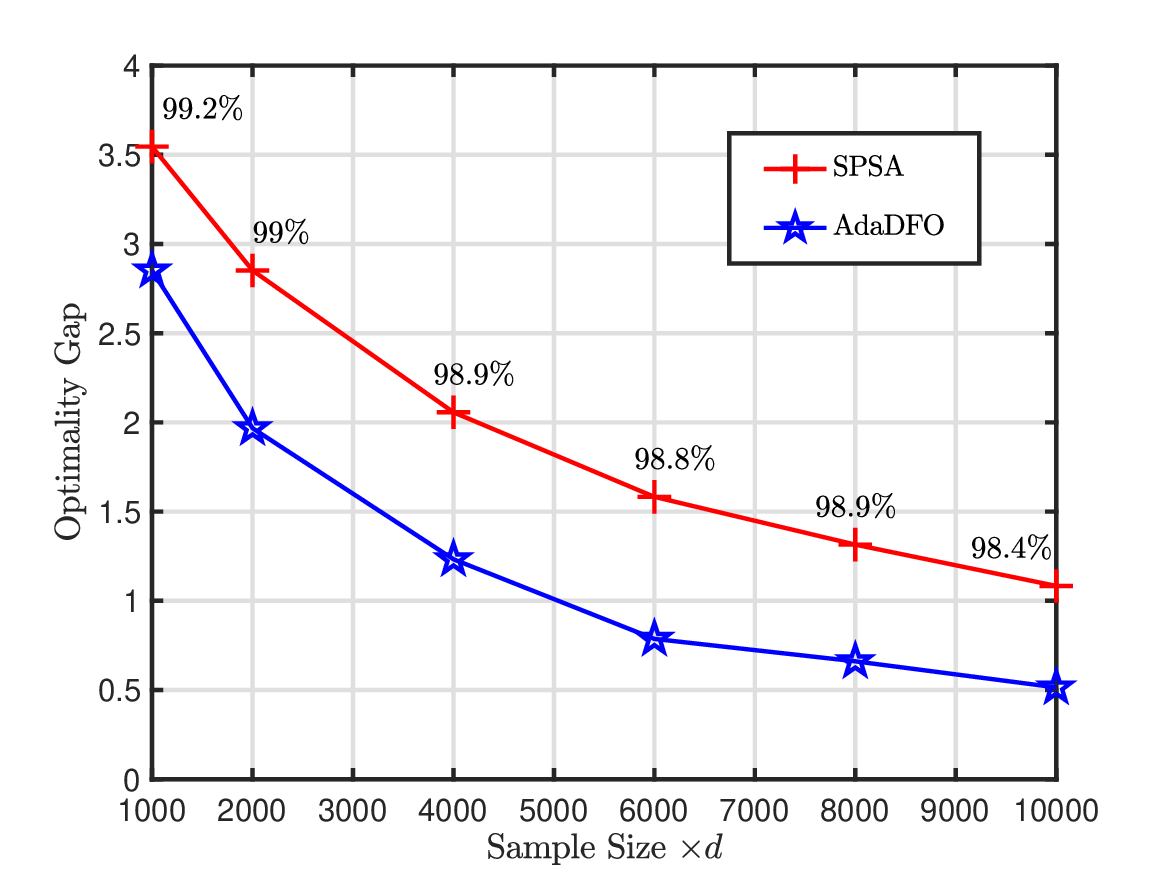}
	\label{fig:OG_Rosenbrock}
\end{figure}

Figure \ref{fig:OG_Rosenbrock} illustrates the average OG between our algorithm and SPSA on the Rosenbrock function. For this type of problem, we set ${\bx}_0 = (-1.9, 2)^{\top}$. Figure \ref{fig:OG_Rosenbrock} includes 6 numbers representing the success rate of the SPSA algorithm across 1000 macroreplications, where the success rate is defined as the probability that the OG is less than $F({\bx}_0) - F({\bx}^*)$. The average OG of SPSA is calculated as the mean of the OGs from all successful macroreplications. As shown in the figure, our algorithm provides better results relative to the SPSA algorithm. Furthermore, our algorithm consistently avoids failures, steadily converges to the optimal solution, and is more robust than the SPSA algorithm.

Similar results are exhibited in Table \ref{tab:Problem3}, where we optimize the third type of problem from the initial point ${\bx}_0 = (3, 1, ..., 3, 1)^{\top}$. Note that the number of sample pairs used is $d$ times the value listed in the table and the initial gap $F({\bx}_0) - F({\bx}^*)$ is approximately $10^{8}$, which is extremely large.
It follows from Table \ref{tab:Problem3} that our algorithm is remarkably successful, giving results that are better than SPSA. For instance, when the sample pairs are $1000d$ and the noise level is 1, the average OG of SPSA is approximately $10^5$ times higher than that of our algorithm. Additionally, SPSA struggles to optimize functions with sharp changes in steepness, and its performance remains almost the same across different noise levels. The reason is that in the initial stage, the step size must be very small to ensure convergence. Otherwise, the algorithm will diverge from the optimal solution \footnote{To clarify, we assume no noise. Then, at ${\bx}_0$, the true gradient is approximately $5.24\times 10^6 \times (1,-1,...,1,-1)^{\top}$. If the step size is too large, the predicted next step will be significantly farther from the optimal solution $(1, 1,...,1)^{\top}$.}. As a result, when the algorithm reaches the flat region, it is difficult to make further progress towards the optimum due to the small step size. In contrast, our algorithm employs a line search that gradually shrinks the step size from $\tilde{a}$ at each iteration, effectively avoiding this issue.

\begin{table}[t!]
\renewcommand\arraystretch{0.9}
\centering
\caption{Comparison of Algorithm \ref{alg:cor_cfd_dfo_ls} with SPSA algorithm for the third type of problem.}
\label{tab:Problem3}
\vspace{1em}
\begin{tabular}{llllllll}
\toprule
 &  & \multicolumn{3}{c}{Solution error} & \multicolumn{3}{c}{Optimality gap} \\
\cmidrule(lr){3-5} \cmidrule(lr){6-8}
$\sigma({\bx})$ & Method & \multicolumn{1}{r}{1,000} & \multicolumn{1}{r}{5,000} & \multicolumn{1}{r}{10,000} & \multicolumn{1}{r}{1,000} & \multicolumn{1}{r}{5,000} & \multicolumn{1}{r}{10,000} \\
\midrule
\multirow{2}{*}{0.1} & AdaDFO & \multicolumn{1}{r}{4.42} & \multicolumn{1}{r}{3.49} & \multicolumn{1}{r}{3.09} & \multicolumn{1}{r}{0.37} & \multicolumn{1}{r}{0.11} & \multicolumn{1}{r}{0.07} \\
& SPSA       & \multicolumn{1}{r}{8.98} & \multicolumn{1}{r}{8.86} & \multicolumn{1}{r}{8.84} & \multicolumn{1}{r}{2.17$\times 10^5$} & \multicolumn{1}{r}{1.12$\times 10^5$} & \multicolumn{1}{r}{8.53$\times 10^4$} \\
\multirow{2}{*}{1}   & AdaDFO & \multicolumn{1}{r}{5.84} & \multicolumn{1}{r}{4.40} & \multicolumn{1}{r}{3.68} & \multicolumn{1}{r}{3.59} & \multicolumn{1}{r}{1.01} & \multicolumn{1}{r}{0.62} \\
& SPSA       & \multicolumn{1}{r}{8.98} & \multicolumn{1}{r}{8.86} & \multicolumn{1}{r}{8.84} & \multicolumn{1}{r}{2.17$\times 10^5$} & \multicolumn{1}{r}{1.11$\times 10^5$} & \multicolumn{1}{r}{8.50$\times 10^4$} \\
\multirow{2}{*}{10}  & AdaDFO & \multicolumn{1}{r}{6.70} & \multicolumn{1}{r}{5.64} & \multicolumn{1}{r}{4.90} & \multicolumn{1}{r}{18.19} & \multicolumn{1}{r}{10.26} & \multicolumn{1}{r}{7.48} \\
& SPSA       & \multicolumn{1}{r}{8.98} & \multicolumn{1}{r}{8.86} & \multicolumn{1}{r}{8.84} & \multicolumn{1}{r}{2.17$\times 10^5$} & \multicolumn{1}{r}{1.11$\times 10^5$} & \multicolumn{1}{r}{8.50$\times 10^4$} \\
\bottomrule
\end{tabular}
\end{table}

\section{Conclusions}\label{sec:Conclusion}
In this paper, we address the limitations of batch-based optimization in gradient-based stochastic search methods for DFO. Specifically, to overcome the challenges associated with selecting the batch size and constructing an effective batch-based gradient estimator, we propose a novel approach that combines the Cor-FD gradient estimate with an adaptive sampling condition. This combination allows us to obtain an appropriate gradient surrogate for KW-type stochastic approximation method.
We prove that, under mild conditions, the use of a properly chosen constant step size ensures convergence. Additionally, we derive the sample complexity of our method, which demonstrates that its convergence rate does not deteriorate compared to the KWSA method. In the black-box scenario, we introduce a new stochastic line search technique to adaptively tune the step size.
Numerical experiments confirm the effectiveness of our proposed algorithm, showing that it outperforms both the KWSA and SPSA algorithms when solving DFO problems.

\newpage
\bibliographystyle{plainnat}
\bibliography{mybibfile}

\begin{thebibliography}{29}
\providecommand{\natexlab}[1]{#1}
\providecommand{\url}[1]{\texttt{#1}}
\expandafter\ifx\csname urlstyle\endcsname\relax
  \providecommand{\doi}[1]{doi: #1}\else
  \providecommand{\doi}{doi: \begingroup \urlstyle{rm}\Url}\fi

\bibitem[Audet and Hare(2017)]{audet2017derivative}
Charles Audet and Warren Hare.
\newblock \emph{Derivative-free and blackbox optimization}.
\newblock Springer, 2017.

\bibitem[Barton and Ivey~Jr(1996)]{barton1996nelder}
Russell~R Barton and John~S Ivey~Jr.
\newblock Nelder-mead simplex modifications for simulation optimization.
\newblock \emph{Management Science}, 42\penalty0 (7):\penalty0 954--973, 1996.

\bibitem[Berahas et~al.(2019)Berahas, Byrd, and Nocedal]{berahas2019derivative}
Albert~S Berahas, Richard~H Byrd, and Jorge Nocedal.
\newblock Derivative-free optimization of noisy functions via quasi-newton
  methods.
\newblock \emph{SIAM Journal on Optimization}, 29\penalty0 (2):\penalty0
  965--993, 2019.

\bibitem[Bollapragada et~al.(2018)Bollapragada, Nocedal, Mudigere, Shi, and
  Tang]{bollapragada2018progressive}
Raghu Bollapragada, Jorge Nocedal, Dheevatsa Mudigere, Hao-Jun Shi, and Ping
  Tak~Peter Tang.
\newblock A progressive batching l-bfgs method for machine learning.
\newblock In \emph{International Conference on Machine Learning}, pages
  620--629. PMLR, 2018.

\bibitem[Bollapragada et~al.(2024)Bollapragada, Karamanli, and
  Wild]{bollapragada2024derivative}
Raghu Bollapragada, Cem Karamanli, and Stefan~M Wild.
\newblock Derivative-free optimization via adaptive sampling strategies.
\newblock \emph{arXiv preprint arXiv:2404.11893}, 2024.

\bibitem[Broadie et~al.(2011)Broadie, Cicek, and Zeevi]{broadie2011general}
Mark Broadie, Deniz Cicek, and Assaf Zeevi.
\newblock General bounds and finite-time improvement for the
  \text{Kiefer-Wolfowitz} stochastic approximation algorithm.
\newblock \emph{Operations Research}, 59\penalty0 (5):\penalty0 1211--1224,
  2011.

\bibitem[Chang et~al.(2013)Chang, Hong, and Wan]{chang2013stochastic}
Kuo-Hao Chang, L~Jeff Hong, and Hong Wan.
\newblock Stochastic trust-region response-surface method (strong)—a new
  response-surface framework for simulation optimization.
\newblock \emph{INFORMS Journal on Computing}, 25\penalty0 (2):\penalty0
  230--243, 2013.

\bibitem[Fazel et~al.(2018)Fazel, Ge, Kakade, and Mesbahi]{fazel2018global}
Maryam Fazel, Rong Ge, Sham Kakade, and Mehran Mesbahi.
\newblock Global convergence of policy gradient methods for the linear
  quadratic regulator.
\newblock In \emph{International Conference on Machine Learning}, pages
  1467--1476. PMLR, 2018.

\bibitem[Fox and Glynn(1989)]{Fox1989Replication}
Bennett~L. Fox and Peter~W. Glynn.
\newblock Replication schemes for limiting expectations.
\newblock \emph{Probability in the Engineering and Informational Sciences},
  3\penalty0 (3):\penalty0 299--318, 1989.

\bibitem[Ghadimi and Lan(2012)]{ghadimi2012optimal}
Saeed Ghadimi and Guanghui Lan.
\newblock Optimal stochastic approximation algorithms for strongly convex
  stochastic composite optimization i: A generic algorithmic framework.
\newblock \emph{SIAM Journal on Optimization}, 22\penalty0 (4):\penalty0
  1469--1492, 2012.

\bibitem[Golovin et~al.(2017)Golovin, Solnik, Moitra, Kochanski, Karro, and
  Sculley]{Golovin2017GoogleVizier}
Daniel Golovin, Benjamin Solnik, Subhodeep Moitra, Greg Kochanski, John Karro,
  and David Sculley.
\newblock Google vizier: A service for black-box optimization.
\newblock In \emph{Proceedings of the 23rd ACM SIGKDD International Conference
  on Knowledge Discovery and Data Mining}, pages 1487--1495, 2017.

\bibitem[Hong and Zhang(2021)]{hong2021surrogate}
L~Jeff Hong and Xiaowei Zhang.
\newblock Surrogate-based simulation optimization.
\newblock In \emph{Tutorials in Operations Research: Emerging Optimization
  Methods and Modeling Techniques with Applications}, pages 287--311. INFORMS,
  2021.

\bibitem[Hu and Fu(2024)]{hu2024convergence}
Jiaqiao Hu and Michael~C Fu.
\newblock On the convergence rate of stochastic approximation for
  gradient-based stochastic optimization.
\newblock \emph{Operations Research}, 2024.

\bibitem[Kiefer and Wolfowitz(1952)]{Kiefer1952Stochastic}
Jack Kiefer and Jacob Wolfowitz.
\newblock Stochastic estimation of the maximum of a regression function.
\newblock \emph{The Annals of Mathematical Statistics}, 23\penalty0
  (3):\penalty0 462--466, 1952.

\bibitem[Kim et~al.(2015)Kim, Pasupathy, and Henderson]{kim2015guide}
Sujin Kim, Raghu Pasupathy, and Shane~G Henderson.
\newblock A guide to sample average approximation.
\newblock \emph{Handbook of simulation optimization}, pages 207--243, 2015.

\bibitem[Larson et~al.(2019)Larson, Menickelly, and Wild]{larsonWild2019DFO}
Jeffrey Larson, Matt Menickelly, and Stefan~M. Wild.
\newblock Derivative-free optimization methods.
\newblock \emph{Acta Numerica}, 28:\penalty0 287--404, 2019.

\bibitem[Li and Lam(2020)]{Li2020Optimally}
Haidong Li and Henry Lam.
\newblock Optimally tuning finite-difference estimators.
\newblock In \emph{2020 Winter Simulation Conference (WSC)}, pages 457--468.
  IEEE, 2020.

\bibitem[Liang et~al.(2024)Liang, Liu, and Zhang]{liang2024cor}
Guo Liang, Guangwu Liu, and Kun Zhang.
\newblock A correlation-induced finite difference estimator, 2024.
\newblock URL \url{https://arxiv.org/abs/2405.05638}.

\bibitem[Mania et~al.(2018)Mania, Guy, and Recht]{mania2018simple}
Horia Mania, Aurelia Guy, and Benjamin Recht.
\newblock Simple random search of static linear policies is competitive for
  reinforcement learning.
\newblock \emph{Advances in neural information processing systems}, 31, 2018.

\bibitem[Robbins and Monro(1951)]{Robbins1951Stochastic}
Herbert Robbins and Sutton Monro.
\newblock A stochastic approximation method.
\newblock \emph{The Annals of Mathematical Statistics}, 22\penalty0
  (3):\penalty0 400--407, 1951.

\bibitem[Scheinberg(2022)]{scheinberg2022finite}
Katya Scheinberg.
\newblock Finite difference gradient approximation: To randomize or not?
\newblock \emph{INFORMS Journal on Computing}, 34\penalty0 (5):\penalty0
  2384--2388, 2022.

\bibitem[Schittkowski(2012)]{schittkowski2012more}
Klaus Schittkowski.
\newblock \emph{More Test Examples for Nonlinear Programming Codes}, volume
  282.
\newblock Springer Science \& Business Media, 2012.

\bibitem[Shashaani et~al.(2018)Shashaani, Hashemi, and Raghu]{Shashaani2018}
Sara Shashaani, Fatemeh~S. Hashemi, and Pasupathy Raghu.
\newblock Astro-df: A class of adaptive sampling trust-region algorithms for
  derivative-free stochastic optimization.
\newblock \emph{Optimization Methods and Software}, 28\penalty0 (4):\penalty0
  3145--3176, 2018.

\bibitem[Shi et~al.(2023)Shi, Xuan, Oztoprak, and Nocedal]{shi2023numerical}
Hao-Jun~Michael Shi, Qiming~Melody Xuan, Figen Oztoprak, and Jorge Nocedal.
\newblock On the numerical performance of finite-difference-based methods for
  derivative-free optimization.
\newblock \emph{Optimization Methods and Software}, 38\penalty0 (2):\penalty0
  289--311, 2023.

\bibitem[Spall(1992)]{Spall1992Multivariate}
James~C. Spall.
\newblock Multivariate stochastic approximation using a simultaneous
  perturbation gradient approximation.
\newblock \emph{IEEE Transactions on Automatic Control}, 37\penalty0
  (3):\penalty0 332--341, 1992.

\bibitem[Spall(1997)]{Spall1997Multivariate}
James~C. Spall.
\newblock A one-measurement form of simultaneous perturbation stochastic
  approximation.
\newblock \emph{Automatica}, 33\penalty0 (1):\penalty0 109--112, 1997.

\bibitem[Spall(2005)]{spall2005introduction}
James~C Spall.
\newblock \emph{Introduction to stochastic search and optimization: estimation,
  simulation, and control}.
\newblock John Wiley \& Sons, 2005.

\bibitem[Wang et~al.(2024)Wang, Liang, Liu, and Zhang]{Wang2024DFO}
Du-Yi Wang, Guo Liang, Guangwu Liu, and Kun Zhang.
\newblock Derivative-free optimization via finite difference approximation: An
  experimental study, 2024.
\newblock URL \url{https://arxiv.org/abs/2411.00112}.

\bibitem[Zazanis and Suri(1993)]{Zazanis1993Convergence}
Michael~A. Zazanis and Rajan Suri.
\newblock Convergence rates of finite-difference sensitivity estimates for
  stochastic systems.
\newblock \emph{Operations Research}, 41\penalty0 (4):\penalty0 694--703, 1993.

\end{thebibliography}

\newpage
\appendix
\renewcommand{\appendixname}{Appendix~\Alph{section}}

\section{Proof of Algorithm \ref{alg:cor_cfd_dfo_cs}}
\subsection{Proof of Theorem \ref{thm:Converge_Iteration}} \label{app:proofTheorem2}
Let $\boldsymbol{\eta}_k := {\bx}_k - {\bx}^*$. We can write \eqref{eq:SA} as 
\begin{align*}
	{\boldsymbol{\eta}}_{k+1} = {\boldsymbol{\eta}}_{k} - a_k \nabla F({\bx}_k) - a_k (g_k({\bx}_k) - \nabla F({\bx}_k)).
\end{align*}
Then we have
\begin{align}
	||{\boldsymbol{\eta}}_{k+1}||^2 =& ||{\boldsymbol{\eta}}_{k}||^2 + a_k^2 ||\nabla F({\bx}_k)||^2 + a_k^2||g_k({\bx}_k) - \nabla F({\bx}_k)||^2 - 2a_k {\boldsymbol{\eta}}_{k}^{\top} \nabla F({\bx}_k)\nonumber\\
	&- 2a_k {\boldsymbol{\eta}}_{k}^{\top}(g_k({\bx}_k) - \nabla F({\bx}_k)) + 2a_k^2 \nabla F({\bx}_k)^{\top} (g_k({\bx}_k) - \nabla F({\bx}_k)).\label{eq:Norm_Square}
\end{align}
According to Assumption \ref{ass:SC} and Taylor expansion, we have
\begin{align*}
	F({\bx}^*) \geq F({\bx}_k) - \nabla F({\bx}_k)^{\top}{\boldsymbol{\eta}}_{k} + \frac{m}{2}||{\boldsymbol{\eta}}_{k}||^2,
\end{align*}
which indicates that
\begin{align}\label{eq:Cross_LB}
	{\boldsymbol{\eta}}_{k}^{\top} \nabla F({\bx}_k) \geq \frac{m}{2}||{\boldsymbol{\eta}}_{k}||^2 + (F({\bx}^*) - F({\bx}_k)).
\end{align}
In addition, according to Assumption \ref{ass:LS} and Taylor expansion, we have
\begin{align*}
	F\left({\bx}_k - \frac{1}{M}\nabla F({\bx}_k)\right) \leq & F({\bx}_k) - \frac{1}{M}||\nabla F({\bx}_k)||^2 + \frac{M}{2}\frac{1}{M^2}||\nabla F({\bx}_k)||^2\\
	= & F({\bx}_k) - \frac{1}{2M}||\nabla F({\bx}_k)||^2,
\end{align*}
which indicates that
\begin{align}\label{eq:Grad_UB}
	||\nabla F({\bx}_k)||^2 \leq 2M \left(F({\bx}_k) - F\left({\bx}_k - \frac{1}{M}\nabla F({\bx}_k)\right)\right) \leq 2M(F({\bx}_k) - F({\bx}^*)).
\end{align}
It follows from \eqref{eq:Norm_Square}, \eqref{eq:Cross_LB} and \eqref{eq:Grad_UB} that
\begin{align*}
	||{\boldsymbol{\eta}}_{k+1}||^2 \leq & (1 - m a_k)||{\boldsymbol{\eta}}_{k}||^2 - (2a_k - 2Ma_k^2)(F({\bx}_k) - F({\bx}^*)) + a_k^2||g_k({\bx}_k) - \nabla F({\bx}_k)||^2\\
	&- 2a_k {\boldsymbol{\eta}}_{k}^{\top}(g_k({\bx}_k) - \nabla F({\bx}_k)) + 2a_k^2 \nabla F({\bx}_k)^{\top} (g_k({\bx}_k) - \nabla F({\bx}_k)).
\end{align*}
Taking conditional expectation on both sides gives
\begin{align}\label{eq:CondiExp}
	\E[||{\boldsymbol{\eta}}_{k+1}||^2 | \cF_k] \leq & (1 - m a_k)||{\boldsymbol{\eta}}_{k}||^2 - (2a_k - 2Ma_k^2)(F({\bx}_k) - F({\bx}^*))\nonumber\\
	& + a_k^2\E[||\boldsymbol{b}_k||^2 | \cF_k] + a_k^2\E[||\boldsymbol{\epsilon}_k||^2 | \cF_k] - 2a_k {\boldsymbol{\eta}}_{k}^{\top}\boldsymbol{b}_k + 2a_k^2 \nabla F({\bx}_k)^{\top} \boldsymbol{b}_k\nonumber\\
	\leq & (1 - m a_k)||{\boldsymbol{\eta}}_{k}||^2 - (2a_k - 2Ma_k^2)(F({\bx}_k) - F({\bx}^*)) + a_k^2\E[||\boldsymbol{b}_k||^2 | \cF_k]\nonumber\\
	& + a_k^2\E[||\boldsymbol{\epsilon}_k||^2 | \cF_k] + 2a_k \left(||{\boldsymbol{\eta}}_{k}||^2 \E[||\boldsymbol{b}_k||^2 | \cF_k]\right)^{1/2} + 2a_k^2 \left(||\nabla F({\bx}_k)||^2 \E[||\boldsymbol{b}_k||^2 | \cF_k]\right)^{1/2}\nonumber\\
	\leq & (1 - m a_k)||{\boldsymbol{\eta}}_{k}||^2 - (2a_k - 2Ma_k^2)(F({\bx}_k) - F({\bx}^*)) + 2\theta^2a_k^2||\nabla F({\bx}_k)||^2 \nonumber\\
	&+ 2\theta a_k \left(||{\boldsymbol{\eta}}_{k}||^2 ||\nabla F({\bx}_k)||^2\right)^{1/2} + 2\theta a_k^2 ||\nabla F({\bx}_k)||^2\nonumber\\
	\leq & (1 - m a_k)||{\boldsymbol{\eta}}_{k}||^2 - (2a_k - 2Ma_k^2)(F({\bx}_k) - F({\bx}^*)) + 2\theta^2a_k^2||\nabla F({\bx}_k)||^2 \nonumber\\
	&+ 2\theta M a_k ||{\boldsymbol{\eta}}_{k}||^2 + 2\theta a_k^2 ||\nabla F({\bx}_k)||^2\nonumber\\
	\leq & (1 - m a_k + 2\theta M a_k)||{\boldsymbol{\eta}}_{k}||^2 - 2a_k(1 - a_k (2\theta^2 + 2\theta + 1)M)(F({\bx}_k) - F({\bx}^*)),
\end{align}
where the second inequality is due to the Cauchy-Schwarz inequality, the third inequality is because we have the threshold condition $\max\{\E[||\boldsymbol{b}_k||^2 | \cF_k], \E[||\boldsymbol{\epsilon}_k||^2 | \cF_k]\} \leq \theta^2 ||\nabla F({\bx}_k)||^2$, the last inequality is due to \eqref{eq:Grad_UB}. From \eqref{eq:Grad_UB}, Assumption \ref{ass:LS} and Taylor expansion, we have 
\begin{align*}
	||\nabla F({\bx}_k)||^2 \leq 2M(F({\bx}_k) - F({\bx}^*)) \leq 2M\left(\nabla F({\bx}^*)^{\top} {\boldsymbol{\eta}}_{k} + \frac{M}{2}||{\boldsymbol{\eta}}_{k}||^2\right) = M^2||{\boldsymbol{\eta}}_{k}||^2,
\end{align*}
which is the reason of the fourth inequality of \eqref{eq:CondiExp}.

Because $0 < a_k = a \leq 1 / ((2\theta^2 + 2\theta + 1) M)$ and $F({\bx^*}) = \underset{{\bx} \in \cX}{\min}F({\bx})$, it follows from \eqref{eq:CondiExp} that for any $k \geq 0$
\begin{align*}
	\E[||{\boldsymbol{\eta}}_{k+1}||^2] \leq \left(1 - (m - 2\theta M)a\right)\E[||{\boldsymbol{\eta}}_{k}||^2].
\end{align*}
Then we have
\begin{align*}
	\E\left[||{\bx}_k - {\bx}^*||^2\right] \leq \left(1 - (m - 2\theta M)a\right)^{k} ||{\bx}_0 - {\bx}^*||^2. 
\end{align*}

The proof is complete.

\subsection{Proof of Corollary \ref{cor:Converge_Iteration}} \label{app:proofCorollary}
Let $\boldsymbol{\eta}_k := {\bx}_k - {\bx}^*$. Note that
\begin{align*}
	F({\bx}_k) - F({\bx}^*) = \nabla F({\bx}^*)^{\top}({\bx}_k - {\bx}^*) + \frac{1}{2}H({\bx}_k^{\dag})||{\bx}_k - {\bx}^*||^2 = \frac{1}{2}H({\bx}_k^{\dag})||{\boldsymbol{\eta}}_{k}||^2,
\end{align*}
where ${\bx}_k^{\dag}$ lies on the line segment between ${\bx}_k$ and ${\bx}^*$, and the second equality is because $\nabla F({\bx}^*) = 0$. Then from Assumptions \ref{ass:LS} and \ref{ass:SC}, we have
\begin{align}\label{eq:OG_Bound}
	\frac{m}{2}||{\boldsymbol{\eta}}_{k}||^2 \leq F({\bx}_k) - F({\bx}^*) \leq \frac{M}{2}||{\boldsymbol{\eta}}_{k}||^2.
\end{align}
Taking \eqref{eq:OG_Bound} back into \eqref{eq:CondiExp} gives
\begin{align*}
	\E[||{\boldsymbol{\eta}}_{k+1}||^2 | \cF_k] \leq & (1 - m a_k + 2\theta M a_k)||{\boldsymbol{\eta}}_{k}||^2 - 2a_k(1 - a_k (2\theta^2 + 2\theta + 1)M)(F({\bx}_k) - F({\bx}^*))\\
	\leq & (1 - m a_k + 2\theta M a_k)||{\boldsymbol{\eta}}_{k}||^2 - m a_k||{\boldsymbol{\eta}}_{k}||^2 + (2\theta^2 + 2\theta + 1)M^2a_k^2||{\boldsymbol{\eta}}_{k}||^2\\
	= & (1 - (2m - 2\theta M) a_k + (2\theta^2 + 2\theta + 1)M^2a_k^2)||{\boldsymbol{\eta}}_{k}||^2\\
	\leq & (1 - (m - \theta M) a_k)||{\boldsymbol{\eta}}_{k}||^2,
\end{align*}
where the last inequality holds when $0 < a_k = a \leq (m - \theta M) / ((2\theta^2 + 2\theta + 1) M^2)$. Therefore,
\begin{align*}
	\E[||{\boldsymbol{\eta}}_{k+1}||^2 | \cF_k] \leq (1 - (m - \theta M) a)||{\boldsymbol{\eta}}_{k}||^2.
\end{align*}
Then, we can get
\begin{align*}
	\E\left[||{\bx}_k - {\bx}^*||^2\right] \leq (1 - (m - \theta M) a)^{k} ||{\bx}_0 - {\bx}^*||^2.
\end{align*}

The proof is complete.

\subsection{Proof of Theorem \ref{thm:Sample_Complex}} \label{app:proofTheorem3}
Let $\boldsymbol{\eta}_k := {\bx}_k - {\bx}^*$. Denote $\mathcal{K}(\epsilon)$ by the number of iteration to get an $\epsilon$-accurate solution. From Theorem \ref{thm:Converge_Iteration}, we should let
\begin{align*}
	\left(1 - (m - 2\theta M)a\right)^{\mathcal{K}(\epsilon)} ||\boldsymbol{\eta}_0||^2 \leq \epsilon,
\end{align*}
which is equivalent to
\begin{align}\label{eq:Iteration_Num}
	\mathcal{K}(\epsilon) \geq \frac{\log\epsilon - \log||\boldsymbol{\eta}_0||^2}{\log\bar{a}},
\end{align}
where we denote $1 - (m - 2\theta M)a$ by $\bar{a}$ for the sake of convenience. Then, $\mathcal{S}(\epsilon) = \sum_{k=0}^{\mathcal{K}(\epsilon) - 1} \mathcal{S}_k = 2d\sum_{k=0}^{\mathcal{K}(\epsilon) - 1} n_k$. From Proposition \ref{pro:Cor-CFD}, $\max\{\E[||\boldsymbol{b}_k||^2 | \cF_k], \E[||\boldsymbol{\epsilon}_k||^2 | \cF_k]\} \leq \mathcal{C}n_k^{-2/3}$, where $\mathcal{C}$ is a constant depending on the third derivative $\nabla^3 F({\bx})$, the function noise $\sigma({\bx})$ and the dimension $d$. To satisfy the adaptive sampling condition, we should let $\mathcal{C}n_k^{-2/3} \leq \theta^2 ||\nabla F({\bx}_k)||^2$. Because the third derivative $\nabla^3 F({\bx})$ and the function noise $\sigma({\bx})$ are bounded, taking expectation about ${\bx}_k$ gives $n_k^{-2/3} \leq \bar{\mathcal{C}}\E[||\nabla F({\bx}_k)||^2]$, where $\bar{\mathcal{C}} = \theta^2/\E[\mathcal{C}]$. Therefore, we get $n_k \geq \bar{\mathcal{C}}^{-3/2} \left(\E[||\nabla F({\bx}_k)||^2]\right)^{-3/2}$.
According to Taylor expansion
\begin{align*}
	\nabla F({\bx}_k) = \nabla F({\bx}^*) + H({\bx}_k^{\dag})\boldsymbol{\eta}_k = H({\bx}_k^{\dag})\boldsymbol{\eta}_k,
\end{align*}
where ${\bx}_k^{\dag}$ lies on the line segment between ${\bx}_k$ and ${\bx}^*$, and the second equality is because $\nabla F({\bx}^*) = 0$. Then we have
\begin{align*}
	\E\left[||\nabla F({\bx}_k)||^2\right] = \E\left[\boldsymbol{\eta}_k^{\top}H^{\top}({\bx}_k^{\dag})H({\bx}_i^{\dag})\boldsymbol{\eta}_k\right] \leq M^2\E\left[||\boldsymbol{\eta}_k||^2\right] \leq M^2\bar{a}^k ||\boldsymbol{\eta}_0||^2.
\end{align*}
Therefore,
\begin{align*}
	n_k^{-2/3} \leq \bar{\mathcal{C}}\E[||\nabla F({\bx}_k)||^2] \leq \bar{\mathcal{C}}M^2\bar{a}^k ||\boldsymbol{\eta}_0||^2,
\end{align*}
which gives that
\begin{align*}
	n_k \geq (\bar{\mathcal{C}}M^2\bar{a}^k ||\boldsymbol{\eta}_0||^2)^{-3/2}.
\end{align*}

Because $\mathcal{S}(\epsilon) = 2d\sum_{k=0}^{\mathcal{K}(\epsilon) - 1} n_k$, we have
\begin{align*}
	\mathcal{S}(\epsilon) \geq & 2d(\bar{\mathcal{C}}M^2 ||\boldsymbol{\eta}_0||^2)^{-3/2}\sum_{k=0}^{\mathcal{K}(\epsilon) - 1}\bar{a}^{-3k/2} \geq 2d(\bar{\mathcal{C}}M^2 ||\boldsymbol{\eta}_0||^2)^{-3/2}\int_{0}^{\mathcal{K}(\epsilon) - 1}\bar{a}^{-3u/2} du\\
	 \geq & 2d(\bar{\mathcal{C}}M^2 ||\boldsymbol{\eta}_0||^2)^{-3/2}\int_{0}^{\frac{\log\epsilon - \log||\boldsymbol{\eta}_0||^2}{\log\bar{a}} - 1}\bar{a}^{-3u/2} du = \mathcal{C}_1 \epsilon^{-3/2} + \mathcal{C}_2,
\end{align*}
where the second inequality is due to the relationship between the summation and integral, and
\begin{align*}
	\mathcal{C}_1 = \frac{4d(\bar{\mathcal{C}}M^2)^{-3/2}}{3\bar{a}^{-3/2}\log(1/\bar{a})},\quad \mathcal{C}_2 = \frac{4d(\bar{\mathcal{C}}M^2||\boldsymbol{\eta}_0||^2)^{-3/2}}{3\log(1/\bar{a})}.
\end{align*}

The proof is complete.

\section{Proof of Algorithm \ref{alg:cor_cfd_dfo_ls}}
\subsection{Proof of Theorem \ref{thm:Converge_LS}} \label{app:proofTheorem4}
Note that
\begin{align}\label{eq:Fun_Bound}
	F({\bx}_{k+1}) \leq & F({\bx}_{k}) - a_k g^{\top}_k({\bx}_k)\nabla F({\bx}_k) + \frac{M}{2}a_k^2||g_k({\bx}_k)||^2\nonumber\\
	=& F({\bx}_{k}) - a_k \left(\nabla F({\bx}_k) + g_k({\bx}_k) - \nabla F({\bx}_k)\right)^{\top}\nabla F({\bx}_k)\nonumber\\
	& + \frac{M}{2}a_k^2||\nabla F({\bx}_k) + g_k({\bx}_k) - \nabla F({\bx}_k)||^2\nonumber\\
	=& F({\bx}_{k}) - a_k\left(1 - \frac{Ma_k}{2}\right)||\nabla F({\bx}_k)||^2 + \frac{Ma_k^2}{2}||g_k({\bx}_k) - \nabla F({\bx}_k)||^2\nonumber\\
	&- a_k(1 - Ma_k)\left(g_k({\bx}_k) - \nabla F({\bx}_k)\right)^{\top}\nabla F({\bx}_k)\nonumber\\
	\leq & F({\bx}_{k}) - a_k\left(1 - \frac{Ma_k}{2}\right)||\nabla F({\bx}_k)||^2 + \frac{Ma_k^2}{2}||g_k({\bx}_k) - \nabla F({\bx}_k)||^2\nonumber\\
	&+ a_k(1 - Ma_k)||g_k({\bx}_k) - \nabla F({\bx}_k)|| ||\nabla F({\bx}_k)||\nonumber\\
	\leq &F({\bx}_{k}) - a_k\left(1 - \frac{Ma_k}{2}\right)||\nabla F({\bx}_k)||^2 + \frac{Ma_k^2}{2}||g_k({\bx}_k) - \nabla F({\bx}_k)||^2\nonumber\\
	&+ \frac{a_k(1 - Ma_k)}{2}\left(||g_k({\bx}_k) - \nabla F({\bx}_k)||^2 + ||\nabla F({\bx}_k)||^2\right)\nonumber\\
	= &F({\bx}_{k}) - \frac{a_k}{2}||\nabla F({\bx}_k)||^2 + \frac{a_k}{2}||g_k({\bx}_k) - \nabla F({\bx}_k)||^2\nonumber\\
	\leq & F({\bx}_{k}) - \frac{a_k(1 - \theta^2)}{2}||\nabla F({\bx}_k)||^2,
\end{align}
where the first inequality is due to Assumption \ref{ass:LS} and the application of Taylor expansion on \eqref{eq:SA}, the second inequality is due to the Cauchy-Schwarz inequality and $a_k \leq 1/M$, the third inequality is because for any two vector ${\bf{u}}, {\bf{v}} \in \R^d$, $||{\bf{u}}|| ||{\bf{v}}|| \leq (||{\bf{u}}||^2 + ||{\bf{v}}||^2)/2$ and the last inequality is because the condition $||g_k({\bx}_k) - \nabla F({\bx}_k)|| \leq \theta |\nabla F({\bx}_k)|$.

Under Assumption \ref{ass:Bound_Error_Fun}, for any ${\bx}$, we have 
\begin{align*}
	f({\bx}) - \bar{\epsilon}_f \leq F({\bx}) \leq f({\bx}) - \bar{\epsilon}_f.
\end{align*}
Therefore, it follows from \eqref{eq:Fun_Bound} that
\begin{align}\label{eq:LS_Verify}
	f({\bx}_{k+1}) \leq f({\bx}_{k}) - \frac{a_k(1 - \theta^2)}{2}||\nabla F({\bx}_k)||^2 + 2\bar{\epsilon}_f \leq f({\bx}_{k}) - \frac{a_k(1 - \theta)}{2(1 + \theta)}||g_k({\bx}_k)||^2 + 2\bar{\epsilon}_f,
\end{align}
where the second inequality is because $||g_k({\bx}_k)|| \leq (1 + \theta) ||\nabla F({\bx}_k)||$. From \eqref{eq:LS_Verify}, we find that the stochastic Armijo condition holds when $a_k \leq 1/M$ and $l_1 \leq (1 - \theta) / (2 (1 + \theta))$. Thus, the stochastic Armijo condition outputs a step size $a_k = \tilde{a}$ or $a_k > l_2/M$. Then we have
\begin{align}\label{eq:LS_Converge}
	f({\bx}_{k+1}) \leq f({\bx}_{k}) - l_1 \min\{\tilde{a}, l_2/M\}||g_k({\bx}_k)||^2 + 2\bar{\epsilon}_f.
\end{align}
From $||g_k({\bx}_k)|| \geq (1 - \theta) ||\nabla F({\bx}_k)||$ and $F({\bx}) - \bar{\epsilon}_f \leq f({\bx}) \leq F({\bx}) - \bar{\epsilon}_f$, it follows from \eqref{eq:LS_Converge} that
\begin{align}\label{eq:LS_Result_Iteration}
	F({\bx}_{k+1}) \leq F({\bx}_{k}) - \min\{l_1\tilde{a}(1-\theta)^2, l_1l_2(1 - \theta)^2/M\}||\nabla F({\bx}_k)||^2 + 4\bar{\epsilon}_f.
\end{align}
In addition, according to Assumption \ref{ass:LS} and Taylor expansion, we have
\begin{align*}
	F\left({\bx}_k - \frac{1}{m}\nabla F({\bx}_k)\right) \geq & F({\bx}_k) - \frac{1}{m}||\nabla F({\bx}_k)||^2 + \frac{m}{2}\frac{1}{m^2}||\nabla F({\bx}_k)||^2\\
	= & F({\bx}_k) - \frac{1}{2m}||\nabla F({\bx}_k)||^2,
\end{align*}
which indicates that
\begin{align}\label{eq:Grad_LB}
	||\nabla F({\bx}_k)||^2 \geq 2m \left(F({\bx}_k) - F\left({\bx}_k - \frac{1}{m}\nabla F({\bx}_k)\right)\right) \geq 2m(F({\bx}_k) - F({\bx}^*)).
\end{align}
Combining \eqref{eq:LS_Result_Iteration} and \eqref{eq:Grad_LB} gives 
\begin{align*}
	F({\bx}_{k+1}) - F({\bx}^*) \leq \left(1 - \min\{2ml_1\tilde{a}(1-\theta)^2, 2ml_1l_2(1 - \theta)^2/M\}\right)(F({\bx}_k) - F({\bx}^*)) + 4\bar{\epsilon}_f.
\end{align*}
Denote $\rho$ by $1 - \min\{2ml_1\tilde{a}(1-\theta)^2, 2ml_1l_2(1 - \theta)^2/M\}$, we have
\begin{align*}
	F({\bx}_{k+1}) - F({\bx}^*) - \frac{4\bar{\epsilon}_f}{1 - \rho} \leq \rho \left(F({\bx}_k) - F({\bx}^*) - \frac{4\bar{\epsilon}_f}{1 - \rho}\right).
\end{align*}
Therefore,
\begin{align*}
	F({\bx}_k) - F({\bx}^*) - \frac{4\bar{\epsilon}_f}{1 - \rho} \leq \rho^{k} \left(F({\bx}_0) - F({\bx}^*) - \frac{4\bar{\epsilon}_f}{1 - \rho}\right).
\end{align*}

The proof is complete.

\section{Parameter Settings in Cor-CFD Method.}
The complete procedure of Cor-CFD method can be found in \cite{liang2024cor}. When applying Cor-CFD in Algorithm \ref{alg:cor_cfd_dfo_ls}, we set the number of perturbations as $K = 5$ and the number of bootstraps as $I = 100$. We set the initial perturbation generator as the truncated normal distribution. At $k$-th iteration, for the first and second types of problems, the mean is 0, variance is $n_k^{-1/5}$, where $n_k$ is the number of sample pairs, and truncated interval is $\left[0.1n_k^{-1/5}, \infty\right)$. For the third type of problem, the mean is 0, variance is $0.1n_k^{-1/5}$ and truncated interval is $\left[0.01n_k^{-1/5}, \infty\right)$.

\end{document}